\documentclass[10pt]{article}

\usepackage[english]{babel}
\usepackage[a4paper,top=20mm,bottom=20mm,left=25mm,right=25mm,marginparwidth=17.5mm]{geometry}

\usepackage{algorithm,algpseudocode}
\usepackage{amsmath}
\usepackage{amssymb}
\usepackage[backend=biber, style=numeric]{biblatex}
\usepackage{caption}
\usepackage{csquotes}
\usepackage{graphicx}
\usepackage[colorlinks=true, allcolors=blue]{hyperref}
\usepackage{listings}
\usepackage[version=4]{mhchem}
\usepackage{microtype}
\usepackage{standalone}
\usepackage{subcaption}
\expandafter\def\csname ver@etex.sty\endcsname{} 
\usepackage{tikz-network}

\usepackage{mycolors}
\usepackage{mycommands}

\usetikzlibrary{positioning}
\begin{document}
\title{A hierarchical dynamical low-rank algorithm for the stochastic description of large reaction networks}
\author{Lukas Einkemmer\thanks{Department of Mathematics, Universit\"{a}t Innsbruck, Innsbruck, Tyrol, Austria} \and Julian Mangott\footnotemark[1]\;\,\footnote{julian.mangott@uibk.ac.at} \and Martina Prugger\thanks{Max-Planck-Institut f\"{u}r Plasmaphysik, Garching, Bavaria, Germany}
}
\maketitle

\begin{abstract}
The stochastic description of chemical reaction networks with the kinetic chemical master equation (CME) is important for studying biological cells, but it suffers from the curse of dimensionality: The amount of data to be stored grows exponentially with the number of chemical species and thus exceeds the capacity of common computational devices for realistic problems. Therefore, time-dependent model order reduction techniques such as the dynamical low-rank approximation are desirable.
In this paper we propose a dynamical low-rank algorithm for the kinetic CME using binary tree tensor networks. The dimensionality of the problem is reduced in this approach by hierarchically dividing the reaction network into partitions. Only reactions that cross partitions are subject to an approximation error.
We demonstrate by two numerical examples (a 5-dimensional lambda phage model and a 20-dimensional reaction cascade) that the proposed method drastically reduces memory consumption and shows improved computational performance and better accuracy compared to a Monte Carlo method.
\end{abstract}

\section{Introduction}
The simulation of chemical reaction networks is an essential tool for gaining insights into (bio)chemical processes. These reaction systems are often described by deterministic rate equations, which amounts to solving a system of ordinary differential equations (ODEs) for the time-dependent concentrations of the chemical species in the system. However, when only few particles are present, the approximation of discrete particle numbers with continuous concentrations is not appropriate and moreover the deterministic ODE description does not model stochastic fluctuations. There are many effects such as stochastic switching, focusing or resonance, which can not be described by deterministic rate equations \cite{Erban_2020}; therefore, stochastic effects are often needed to describe important features of realistic systems in biochemistry \cite{Tonn_2019,Grima_2008,Niepel_2009,Paszek_2010}.

The fundamental model of the stochastic description of chemical reaction systems is the chemical master equation (CME). The CME describes the time evolution of a probability density function that, in addition to time, depends on the population numbers. Several analytical results for stationary distributions are known (e.g. \cite{VanKampen_1976,Anderson_2010}), but for the time-dependent case exact solutions only exist for monomolecular reactions \cite{Jahnke_2007} and for hierarchic first-order networks (i.e. monomolecular reactions, (auto-)catalytic and splitting reactions) \cite{Reis_2018}.

In most realistic situations the CME has to be solved numerically. Unfortunately, the CME suffers from the so-called curse of dimensionality; that is, the memory requirement and computational cost scales exponentially in the number of dimensions. Thus solving the CME numerically with a direct method is extremely expensive.
Besides effective models (such as replacing the CME by a chemical Langevin equation \cite{Gillespie_2000} or by introducing rate equations for averaged population numbers \cite{Chen_2010}), essentially two classes of solution methods can be discerned in the literature to circumvent the curse of dimensionality: The first class are Monte Carlo methods, whereas the second class comprehends complexity reduction and sparse methods.

The de facto standard in the class of Monte Carlo methods is the stochastic simulation algorithm (SSA), which was first described in \cite{Gillespie_1976,Gillespie_1977}. In this approach trajectories are sampled in accordance with the solution of the CME (the yet unknown probability distribution). This method has been improved over the course of years, for example by incorporating several reactions into one sampling period. This so-called $\tau$-leaping method was first proposed in \cite{Gillespie_2001} and is accurate when the reaction rates do not change much during a time step. Improvements of this approximation have been made for example by deriving an upper bound for the number of reactions in one sampling period \cite{Albi_2022}.
Since SSA is a Monte Carlo approach, it overcomes the curse of dimensionality by sampling the high dimensional state space. However, in order to reconstruct the underlying probability distribution, typically many trajectories have to be calculated. This is the main downside of this approach: it only converges slowly (as $1/\sqrt{N}$, where $N$ is the number of trajectories) and, in particular, fine details of the probability distribution are not well resolved when not enough samples are used. However, in many systems interesting events only occur rarely \cite{Kazeev_2014} and producing many trajectories is computationally expensive. Additionally, a post-processing step is typically required for obtaining the probability distribution or observables (such as moments of the distribution) from the individual trajectories.

The methods of the second class aim to solve the CME directly through complexity reduction and sparse matrix techniques. There is a cornucopia of different algorithms described in the literature and we thus can only give a non-exhaustive overview. Since the CME is an infinite systems of ODEs (there are infinitely many possible population numbers), it is necessary to reduce it to a finite state space and thus all values of the probability distribution are truncated when they are zero up to a prescribed tolerance. This approach is called finite-state projection (FSP) method and was first described in \cite{Munsky_2006}. An improvement of this technique is the use of sliding windows, where the finite state space follows the dynamics of the probability function \cite{Henzinger_2009,Wolf_2010,Dinh_2020}. Unfortunately, these methods do not reduce the dimensionality of the problem, since the truncation values of the population number for every dimension may be large and the exponential scaling with the number of dimensions for memory and computational cost is still present.
Therefore, many other methods have been built on top of FSP, such as sparse grids \cite{Hegland_2007,Hegland_2008}, Wavelet compression \cite{Jahnke_2010,Jahnke_2010a}, spectral methods \cite{Engblom_2008} or Krylov subspace methods \cite{MacNamara_2008}.

In this paper we apply a low-rank approximation to overcome the curse of dimensionality. The main idea of the low-rank approximation is the description of the problem by a small set of lower-dimensional basis functions (the so-called low-rank factors), which are combined to approximate the high-dimensional probability distribution. 
Since we are using only a small number of basis functions (the so-called rank of the approximation), the memory requirements and the computational effort are reduced significantly. 
In order to advance the approximation forward in time, evolution equations are derived for the individual degrees of freedom of the low-rank approximation.
More specifically, the derivative of the approximation is projected onto the tangent space of the low-rank manifold, and thus ensures that the approximation itself always remains low-rank \cite{Lubich_2008}. The method proposed in \cite{Lubich_2008} suffers unfortunately from instability when the singular values of the approximation are very small. Small singular values are however required for a good approximation. Fortunately, this problem was overcome by the projector-splitting \cite{Lubich_2014} and later the Basis Update \& Galerkin (BUG) integrators \cite{Ceruti_2022}, which are robust with respect to small singular values. Both integrators have been improved subsequently by proposing schemes which are second-order \cite{Ceruti_2024b}, asymptotic-preserving and conservative \cite{Einkemmer_2024,Einkemmer_2022,Einkemmer_2021b}, parallel \cite{Kusch_2024,Ceruti_2024a} or rank-adaptive \cite{Ceruti_2022a}. In particular, the projector-splitting and BUG integrators were also proposed for Tucker tensors \cite{Lubich_2018} and tree tensor networks (TTNs) \cite{Ceruti_2020,Ceruti_2023}.

We are extending our dynamical low-rank matrix integrator for the kinetic CME, as described in \cite{Einkemmer_2018}, to a more general binary TTN integrator. In \cite{Einkemmer_2018} the main idea was to separate the physical dimensions (i.e., the species) of the problem. The reaction network is split into two partitions, and every partition is then described by a small number of low-rank factors. These low-rank factors are lower-dimensional, as they only depend on the species in their partition. It is thereby guaranteed that all reactions inside of a partition are treated exactly and an approximation is only performed if a reaction involves species of both partitions, i.e., when the reaction crosses the partition boundary. This approach of separating physical dimensions has been used for Boolean models in biology \cite{Prugger_2023}, for problems in plasma physics (see, e.g., \cite{Einkemmer_2018,Cassini_2021,Coughlin_2022,Einkemmer_2020}), radiation transport (see, e.g., \cite{Peng_2020,Peng_2021,Kusch_2021,Einkemmer_2021,Einkemmer_2021a,Kusch_2022}) and quantum spin systems \cite{Ceruti_2024,Sulz_2024}. In these publications the partitioning was done on the level of physical dimensions (in contrast to quantized tensor train (QTT) methods, which we describe below) and the partitions were chosen such that they suit the physical problem (either with a decomposition into spatial and velocity scales, as in \cite{Peng_2020,Kusch_2021,Einkemmer_2020,Einkemmer_2018}, or in coordinates parallel and perpendicular to the magnetic field, as in \cite{Einkemmer_2023}).

The solutions obtained by the dynamical low-rank matrix integrator for the kinetic CME are, in contrast to Monte Carlo methods, completely noise free and are stored in a compressed format. Observables such as moments of the probability distribution can be easily computed from this compressed representation. Another benefit is the ability to keep species which are strongly correlated together without introducing any error. The computational and memory savings of this approach come form the fact that the full solution is approximated by using only a small number of lower-dimensional low-rank factors. We emphasize that the rank (i.e. the number) of the of low-rank factors and their dimension strongly depends on the specific problem and on the partitioning of the reaction network.

Our previous approach of separating the reaction network into two partitions alleviates the curse of dimensionality, but for many dimensions the computational effort and memory requirements may still be prohibitively large, since both partitions still contain at best $d/2$ different species, where $d$ is the total number of species. It is therefore often more preferable to separate the two partitions again in smaller subpartitions. In principle, this partitioning process could be done recursively until reaching the level where each partition contains only a single species (i.e., a single physical dimension). This new method gives now complete flexibility in separating species or keeping them together, depending on the specific problem.
Similar to our previous work, the low-rank factors for each (sub-)partition only depend on the species in their partition and thus are again lower-dimensional. They will be updated with the projector-splitting integrator for tree tensor networks \cite{Ceruti_2020}.

In the past, there have been several other attempts of solving the kinetic CME with a low-rank approximation. In \cite{Jahnke_2008,Hegland_2011} the low-rank factors are allowed to depend only on a single species. This is reminiscent of approximations in quantum mechanics, where single-orbital basis functions are used, which only depend on the coordinates of a single electron, and then are combined to obtain an approximation to the high-dimensional wave function. Another very popular tool in quantum mechanics, the quantized tensor trains (QTTs), have also been applied to the CME \cite{Kazeev_2014,Kazeev_2015,Dolgov_2015,Dinh_2020,Ion_2021}. Like our approach, these QTT methods decompose the reaction network into smaller partitions, but the crucial difference is that not only the physical dimensions of the species are separated, but the physical dimensions are decomposed into even smaller virtual dimensions. Although this reduces the memory requirements further, it comes with several downsides: Species which are coupled by important reaction pathways are inevitably separated and moreover, the propensity functions on the right-hand side of the CME have to be approximated. Given the intricate structures of complex biological networks \cite{Barabasi_2009}, it is doubtful that in biological applications we can consider each species independently, not to mention decomposing single species into virtual dimensions, while still obtaining an accurate approximation with a small rank. Evidence of this gives a numerical example in \cite{Kazeev_2014}, where the required rank is indeed very high.

The remainder of this paper is structured as follows. In Section~\ref{sec:cme} we explain the stochastic formulation of chemical reaction networks and the chemical master equation and we establish our notation. We then recapitulate the main ideas for the matrix case found in \cite{Einkemmer_2024} in Section~\ref{sec:matrix}. Next, the matrix integrator will be generalized to tree tensor networks in Section~\ref{sec:ttn}. In Section~\ref{sec:experiments} we investigate the accuracy and efficiency of the method for a number of examples. Finally, we conclude in Section~\ref{sec:outlook}.

\section{Chemical master equation\label{sec:cme}}
We consider a chemical reaction system with $d$ species $S_{0},\ldots,S_{d-1}$. In the stochastic description the population number (i.e. number of molecules) of the $i$-th species at time $t$ is given by a random variable $\mathcal{X}_{i}(t)$ on the discrete state space $\mathbb{N}_{0}$. In the following we denote by $\vb{x}=(x_0,\ldots,x_{d-1})\in\mathbb{N}_{0}^d$ a realization of the random variable $\vb{\mathcal{X}}(t)=\left(\mathcal{X}_{0}(t),\ldots,\mathcal{X}_{d-1}(t)\right)$ at time $t$. The species can interact through $M$ reaction channels $R_{0},\ldots,R_{M-1}$, where every reaction $R_\mu$ is described by two quantities: The first quantity is the stoichiometric vector $\vb{\nu}_{\mu}=(\nu_{\mu,0},\ldots,\nu_{\mu,d-1})\in\mathbb{Z}^{N}$, which describes the change of the population number $\vb{x}$ caused by reaction $R_\mu$. The second quantity is the propensity function $\alpha_{\mu}(\vb{x}):\mathbb{N}_{0}^d\to[0,\,\infty)$, where $\alpha_\mu(\vb{x})\diff{t}$ is the probability that the reaction $R_\mu$ occurs in the infinitesimally small time interval $[t, t+\diff{t})$.

We will restrict ourselves to well-stirred chemical reaction systems, thus the propensity functions do not show a spatial dependency. Since $\vb{\mathcal{X}}$ is a random variable, we describe the reaction system by the probability density
\begin{equation*}
    P(t,\vb{x})=\mathbb{P}(\mathcal{X}_0(t)=x_0,\ldots,\mathcal{X}_{d-1}(t)=x_{d-1}),
\end{equation*}
which is the solution of the kinetic chemical master equation (CME)
\begin{equation}\label{eq:CME}
    \partial_{t}P(t,\vb{x})=\sum_{\mu=0}^{M-1}\left(\alpha_{\mu}(\vb{x}-\vb{\nu}_{\mu})P(t,\vb{x}-\vb{\nu}_{\mu})-\alpha_{\mu}(\vb{x})P(t,\vb{x})\right).
\end{equation}
The CME can be considered as a infinite system of ODEs with one equation for every state $\vb{x}$. The term ``kinetic'' indicates that the population number can be of any natural number (including 0), in contrast to the special case where only Boolean states are allowed (i.e.~where a species can be either ``activated'' or ``not activated''; see, e.g., \cite{Clarke_2020,Zanudo_2018,Yachie-Kinoshita_2018,Prugger_2023}).

By defining the linear operator
\begin{equation}
    \left(\mathcal{A}P(t,\cdot)\right)(\vb{x})=\sum_{\mu=0}^{M-1}\left(\alpha_{\mu}(\vb{x}-\vb{\nu}_{\mu})P(t,\vb{x}-\vb{\nu}_{\mu})-\alpha_{\mu}(\vb{x})P(t,\vb{x})\right),\label{eq:linear_operator}
\end{equation}
we can write the CME concisely as 
\begin{equation*}
    \partial_{t}P(t,\cdot)=\mathcal{A}P(t,\cdot).
\end{equation*}
Reactions that reduce the population number to negative values are unphysical, therefore we have to omit the arguments $\vb{x}-\vb{\nu}_{\mu}$ in the first term of the right-hand side of Equation~(\ref{eq:CME}) when they become negative. For more details on the stochastic description of reaction systems and the CME in general we refer the reader to, e.g., \cite{Gillespie_1976,Gillespie_1992,Gardiner_2004,Erban_2020}.

The probability distribution satisfies conservation of mass:
\begin{equation}
    \sum_{\vb{x}\in\mathbb{N}_{0}^d}P(t,\vb{x})=\sum_{\vb{x}\in\mathbb{N}_{0}^d}P(0,\vb{x})=1.\label{eq:mass}
\end{equation}
This relation can be directly derived from the CME by integrating Equation~(\ref{eq:CME}) over time and performing a summation over $\vb{x}$, which yields 
\begin{equation}
    \sum_{\vb{x}\in\mathbb{N}_{0}^d}\left(P(t,\vb{x})-P(0,\vb{x})\right)=\text{\ensuremath{\int_{0}^{t}}}\sum_{\mu=0}^{M-1}\sum_{\vb{x}\in\mathbb{N}_{0}^d}\left(\alpha_{\mu}(\vb{x}-\vb{\nu}_{\mu})P(\tilde{t},\vb{x}-\vb{\nu}_{\mu})-\alpha_{\mu}(\vb{x})P(\tilde{t},\vb{x})\right)\,\mathrm{d}\tilde{t}.\label{eq:mass_int}
\end{equation}
The right-hand side of this equation vanishes since
\begin{equation}
    \sum_{\mu=0}^{M-1}\sum_{\vb{x}\in\mathbb{N}_{0}^d}\alpha_{\mu}(\vb{x}-\vb{\nu}_{\mu})P(t,\vb{x}-\vb{\nu}_{\mu})=\sum_{\mu=0}^{M-1}\sum_{\vb{x}\in\mathbb{N}_{0}^d}\alpha_{\mu}(\vb{x})P(t,\vb{x}),\label{eq:substitution}
\end{equation}
which implies the desired result.

Let us illustrate the concepts by a simple example first proposed by Schl\"{o}gl \cite{Schloegl_1971}. The chemical species $S$ in this model is subject to two reversible reactions,
\begin{equation*}
    \ce{3$S$ <=>[$k_0$][$k_1$] 2$S$}, \qquad{} \ce{$S$ <=>[$k_2$][$k_3$] $\star$,}\\
\end{equation*}
where $\star$ denotes chemical species that are of no further interest. Note that in this model the particle number is not conserved, whereas the conservation of probability is still guaranteed. We decompose the two reversible reactions into four unidirectional reactions and assume for this example that all reactions are elementary, i.e., the propensity functions are given by the rate constants times the number of distinct reactant combinations. Let us address the first forward reaction \ce{3$S$ ->[$k_0$] 2$S$} explicitly: Here the number of different combinations of reactants is $\binom{x}{3} = \frac{x (x-1) (x-2)}{6}$. It is common to absorb the factor $6$ into the rate constant $k_0$ and thus the propensity function for this reaction is $\alpha_0(x) = k_0 x(x-1)(x-2)$. Since three particles of $S$ decay into two particles of the same type, the stoichiometric vector for this reaction is $\nu_0 = -1$.
The CME for the four unidirectional, elementary reactions is 
\begin{equation}\label{eq:cme_schloegl}
    \begin{aligned}
        \partial_{t}P(t,x) &= (\alpha_0(x+1)+\alpha_2(x+1)) P(t,x+1) + (\alpha_1(x-1) + \alpha_3(x-1)) P(t,x-1) \\ &\quad - (\alpha_0(x)+\alpha_1(x)+\alpha_2(x)+\alpha_3(x)) P(t,x),
    \end{aligned}
\end{equation}
with propensity functions $\alpha_1(x) = k_1 x(x-1)$, $\alpha_2(x) = k_2 x$ and $\alpha_3(x) = k_3$ for the three remaining reactions. 
It is instructive to compare the stochastic methods for this model with deterministic rate equations: In the deterministic setting we solve the ordinary differential equation
\begin{equation}\label{eq:schloegl_deterministic}
    \frac{\mathrm{d}x}{\mathrm{d}t} = -k_0 x^3 + k_1 x^2 - k_2 x + k_3,
\end{equation}
with parameters $k_0 = 2.5 \cdot 10^{-4}$, $k_1 = 0.18$, $k_2 = 37.5$ and $k_3 = 2200$ taken from \cite{Erban_2020}. With these parameters the rate equation exhibits two stable and one unstable steady state solutions, and starting with an initial condition below or above the separatrix at $x(t)=220$ will end either in the lower or the upper steady state, respectively.
In the stochastic formulation we draw one sample $x(t) \sim P(t, x)$ with the stochastic simulation algorithm (SSA) \cite{Gillespie_1976}.

\begin{figure}
    \centering
    \includegraphics[scale=0.8]{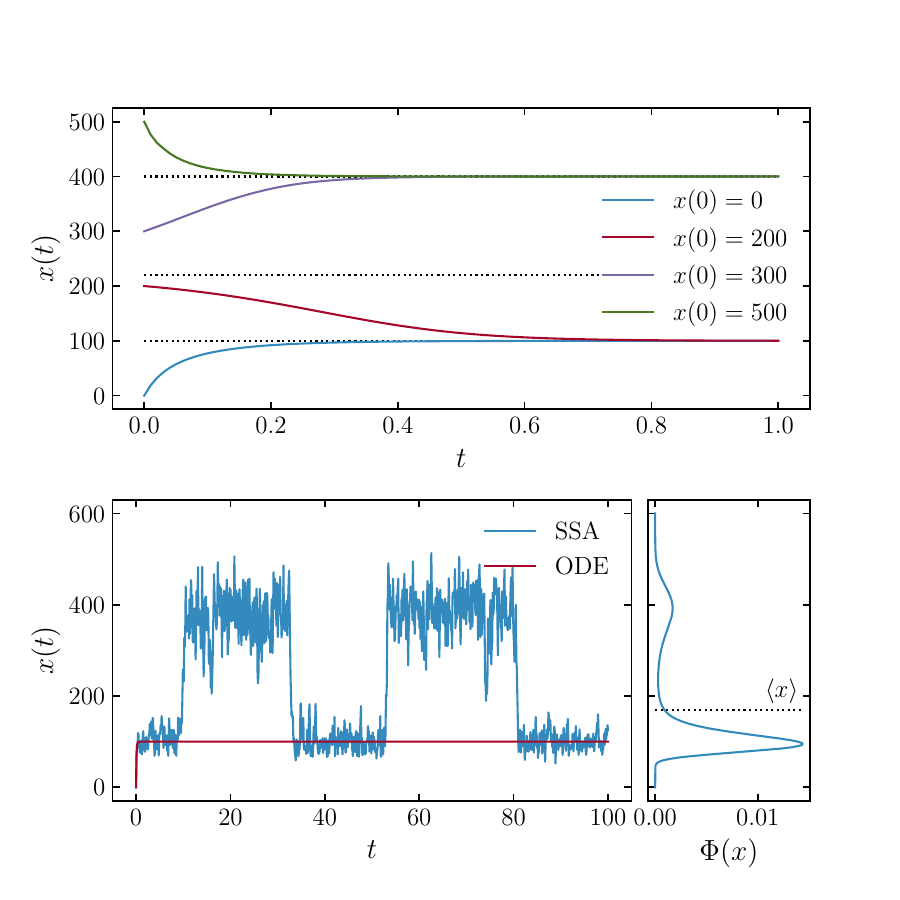}
    \caption{Population number depending on time for the Schl\"{o}gl model with $k_0 = 2.5 \cdot 10^{-4}$, $k_1 = 0.18$, $k_2 = 37.5$ and $k_3 = 2200$ in the deterministic (top) and stochastic setting (bottom, left), with the steady state distribution $\Phi(x) = \lim_{t\to \infty}P(t, x)$ (bottom, right). The deterministic solution was calculated from Equation~(\ref{eq:schloegl_deterministic}) and the single sample $x(t) \sim P(t, x)$ with SSA. Parameters were taken from \cite{Erban_2020}. Note that the steady ensemble average is $\lim_{t \to \infty} \langle x(t) \rangle \approx 169.46$ for the stochastic approach, whereas the rate equation would predict, depending on the initial condition, either $\lim_{t \to \infty} x(t) = 100$ or $400$.\label{fig:schloegl}}
\end{figure}

The results for both methods are shown in Figure~\ref{fig:schloegl}. In the top panel we show the time-dependent population number for four different initial conditions and the three stable steady state solutions are indicated by dashed lines. As expected, the four solutions converge to the two stable steady state solutions.
The plot on the bottom left shows a single trajectory generated by SSA in comparison with the deterministic solution, both for initial condition $x(t)=0$. In the stochastic setting the probability for a transition to the higher steady state is non-zero and therefore the SSA solution fluctuates between the two steady state solutions. This stochastic switching between the stationary states is not captured by the deterministic description, the ODE solution converges to the lower steady state.
The plot on the bottom right shows the steady distribution $\Phi(x) = \lim_{t\to \infty}P(t, x)$. Note that the steady ensemble average is $\lim_{t \to \infty} \langle x(t) \rangle \approx 169.46$ for the stochastic approach with initial condition $P(0,x)=\delta_{x,0}$, whereas the rate equation would predict $\lim_{t \to \infty} x(t) = 100$ with initial condition $x(0)=0$.

We observe that the deterministic approach does not cover stochastic effects, which are however important to model systems in biochemistry accurately \cite{Tonn_2019,Grima_2008,Niepel_2009,Paszek_2010}. We also see from this example, that a single trajectory generated by SSA is by far not enough to reconstruct the full probability distribution. Indeed, SSA only converges with $1/\sqrt{N}$, where $N$ is the number of trajectories, and many trajectories have to be generated to resolve the fine details of the probability distribution. If one is interested in rare events, a direct solution of the CME using complexity reduction or sparse techniques is therefore preferable as we will see.

\section{Recap: Matrix integrator}\label{sec:matrix}
In this section, we recapitulate the main features of the matrix integrator for the kinetic CME \cite{Einkemmer_2024}. We will generalize these ideas later for the tree tensor network case.

As already mentioned in the introduction, the CME can be regarded as an infinite system of ODEs. Thus, we have to truncate the state space to a finite domain in order to turn the CME into a finite problem. This truncation is called finite state projection (FSP) and was first introduced in \cite{Munsky_2006}. The finite domain has to be chosen such that it allows a numerical solution and still captures enough of the information of the full (infinite) system.

The truncated state space is defined as $\Omega_{\vb{\zeta},\vb{\eta}}=\{\vb{x} \in\mathbb{N}_{0}^d:\zeta_{i}\le x_{i}\le\eta_{i}\ \mathrm{for}\ i=0,\dots,d-1\}$, where $\zeta_{i}\in\mathbb{N}_{0}$ and $\eta_{i}\in\mathbb{N}_{0}$ and $\zeta_{i}<\eta_{i}$ ($i=0,\ldots,d-1)$. We denote by $\mathcal{A}_{\vb{\zeta},\vb{\eta}}$ the restriction of the linear operator $\mathcal{A}$ (defined in Equation~\ref{eq:linear_operator}) to $\Omega_{\vb{\zeta},\vb{\eta}}$ and by $P_{\vb{\zeta},\vb{\eta}}(t)$ the solution of the restricted CME $\partial_{t}P_{\vb{\zeta},\vb{\eta}}(t)=\mathcal{A}_{\vb{\zeta},\vb{\eta}}P_{\vb{\zeta},\vb{\eta}}(t)$ with initial condition $P_{\vb{\zeta},\vb{\eta}}(0)$, which is the initial probability distribution restricted to the truncated state space. Defining the total mass $m_{\vb{\zeta},\vb{\eta}}=\sum_{\vb{x}\in\Omega_{\vb{\zeta},\vb{\eta}}}P_{\vb{\zeta},\vb{\eta}}(t,\vb{x})$
and assuming that $m_{\vb{\zeta},\vb{\eta}} \ge 1-\epsilon$, it was shown in \cite{Munsky_2006} that
\begin{equation*}
    P(t,\vb{x})-\epsilon \le P_{\vb{\zeta},\vb{\eta}}(t,\vb{x})\le P(t,\vb{x}) \quad \text{for} \quad \vb{x} \in \Omega_{\vb{\zeta},\vb{\eta}}.
\end{equation*}
These inequalities give an estimation of how close the truncated state space solution approximates the true solution. The main issue of the truncation is how to determine suitable $\vb{\zeta}$ and $\vb{\eta}$ for a given final time $t$ and tolerance $\epsilon>0$, such that $m_{\vb{\zeta},\vb{\eta}} \ge 1-\epsilon$. Solving the reaction network deterministically with ODEs (which is cheap) or biological insight into the system might give a good idea on how to choose $\vb{\zeta}$ and $\vb{\eta}$ a priori. Alternatively, one can implement a scheme with an adaptive truncated state space as proposed in \cite{Henzinger_2009,Wolf_2010,Dinh_2020}.

In what follows, we will always work with the restriction of $P$ on the truncated state space and we omit the truncation values in our notation, i.e., $P \equiv P_{\vb{\zeta},\vb{\eta}}$. Since the memory requirement for storing the truncated probability distribution still scales exponentially with the number of species,
we have to reduce the system size in order to solve the CME using currently available hardware. In \cite{Einkemmer_2024} we have done this by using a dynamical low-rank (DLR) approximation. The main idea of the DLR approximation is to split the species (and thus the reaction network) into two partitions, i.e., $\vb{x} = (\vb{x}^0, \vb{x}^1)$, with $\vb{x}^0 \in \Omega^0_{\vb{\zeta},\vb{\eta}} \subset \mathbb{N}^{d^0}_0$ and $\vb{x}^1 \in \Omega^1_{\vb{\zeta},\vb{\eta}} \subset \mathbb{N}^{d^1}_0$, where $d^0 + d^1 = d$. (Note that in the following, we will denote all quantities belonging to partition $0$ or partition $1$ with superscripts.) Since we work with truncated state spaces (i.e., with a finite number of states), the probability distribution can be represented for a given time $t$ as a $n^0 \times n^1$~matrix, $P(t,\vb{x}) \equiv P(t,\vb{x}^0,\vb{x}^1)$, where $n^0$ and $n^1$ is the size of the truncated state space $\Omega^0_{\vb{\zeta},\vb{\eta}}$ and $\Omega^1_{\vb{\zeta},\vb{\eta}}$, respectively. We are then approximating the probability distribution by $\tilde{r}$ time-dependent low-rank factors $X^0_i(t, \vb{x}^0)$ and $X^1_i(t, \vb{x}^1)$ and a time-dependent, invertible coefficient matrix $S_{ij}(t) \in \mathbb{R}$, which is reminiscent of a truncated singular value decomposition (SVD):
\begin{equation}\label{eq:DLR-approximation}
    P(t, \vb{x}) \approx \sum_{i,j=0}^{\tilde{r}-1} X^0_i(t, \vb{x}^0) S_{ij}(t) X^1_j(t, \vb{x}^1).
\end{equation}
The low-rank factors obey the orthogonality and gauge conditions for the low-rank factors,
\begin{align}
    \langle X_{i}^{0},X_{j}^{0}\rangle_{0}=\delta_{ij} & \quad\textrm{and}\quad\langle X_{i}^{1},X_{j}^{1}\rangle_{1}=\delta_{ij},\qquad\textrm{(orthogonality)},\label{eq:ortho}\\
    \langle X_{i}^{0},\partial_{t}X_{j}^{0}\rangle_{0}=0 & \quad\textrm{and}\quad\langle X_{i}^{1},\partial_{t}X_{j}^{1}\rangle_{1}=0,\qquad\textrm{(gauge condition)},\label{eq:gauge}
\end{align}
where $\delta_{ij}$ denotes the Kronecker delta and $\langle\cdot,\cdot\rangle_{k}$
the inner product on $\ell^{2}(\Omega^k_{\vb{\zeta},\vb{\eta}})$ $(k=0,1)$. In the following, we will omit summation bounds in order to simplify the notation. We agree that sums with latin indices run over all low-rank factors, i.e. $\sum_i\equiv\sum_{i=0}^{\tilde{r}-1}$, and sums with a greek index $\mu$ run over all reactions, i.e. $\sum_\mu\equiv\sum_{\mu=0}^{R-1}$.

Since low-rank factors only depend on the species in one partition (either $\vb{x}^0$ or $\vb{x}^1$), it is ensured that reaction pathways lying within a partition are treated exactly, while reaction pathways that cross the two partitions are taken into account in an approximate way.
The total number of low-rank factors $\tilde{r}$ is called the rank. When the rank is small compared to the partition sizes, i.e., $\tilde{r} \ll n^0$ and $\tilde{r} \ll n^1$, then the total memory requirement is reduced for a given time $t$ from $\bigO(n^0 n^1)$ to $\bigO(\tilde{r}(n^0+n^1)+\tilde{r}^2)$. In the numerical implementation the low-rank factors can be represented as matrices without any further loss of accuracy, since we are working in the truncated state space and the population numbers are discrete. Inner products can be represented by matrix multiplications.

In order to retain efficiency, the full probability distribution must not be formed explicitly by Equation~(\ref{eq:DLR-approximation}). To this end, we have to replace the time evolution of $P(t,\vb{x})$ via the CME by evolution equations for the low-rank factors and the coefficient matrix. In \cite{Einkemmer_2024} we derived schemes with the projector-splitting integrator of \cite{Lubich_2014}, which were first- and second-order accurate in time. Here we will focus on the first-order scheme, for more details on the extension to second-order and on the derivation of the evolution equations we refer to \cite{Einkemmer_2024}.

The initial conditions $X^0_{0,i}(\vb{x}^0)$ and $X^1_{0,j}(\vb{x}^1)$ for the low-rank factors and for the coefficient matrix $S_{0,ij}$ are given by the approximation $P(0,\vb{x}) \approx \sum_{i,j} X^0_{0,i}(\vb{x}^0) S_{0,ij} X^1_{0,j}(\vb{x}^1)$ for the initial value of $P(t,\vb{x})$. When $P(0,\vb{x})$ is in a low-rank form, then obtaining the low-rank factors and the coefficient matrix is straightforward. Otherwise the low-rank factors and the coefficient matrix can be computed from $P(0,\vb{x})$ by a truncated (randomized) SVD.

Our aim now is to obtain a representation for the probability distribution after one time step $\Delta t$. We achieve this by advancing $X^0_{0,i}(\vb{x}^0)$, $S^0_{0,ij}$ and $X^0_{1,j}(\vb{x}^0)$ sequentially in time. 
First, we start by forming $K_i(0,\vb{x}^0)=\sum_j X_{0,j}^0(\vb{x}^0) S_{0,ji}$ and use this as an initial condition for the evolution equation
\begin{equation}\label{eq:K}
    \partial_{t}K_{i}(t,\vb{x}^0)=\sum_\mu\sum_j\left(c_{\mu,ij}^{0}(\vb{x}^0)K_{j}(t,\vb{x}^0-\vb{\nu}^0_{\mu})-d_{\mu,ij}^{0}(\vb{x}^0)K_{j}(t,\vb{x}^0)\right),
\end{equation}
with the time-independent coefficients 
\begin{equation}\label{eq:K_coeff}
    \begin{aligned}
        c_{\mu,ij}^{0}(\vb{x}^0) & =\langle X_{0,i}^{1}(\vb{x}^1),\alpha_{\mu}(\vb{x}^0-\vb{\nu}^0_{\mu},\vb{x}^1-\vb{\nu}^1_{\mu})X_{0,j}^1(\vb{x}^1-\vb{\nu}^1_{\mu})\rangle_{1},\\
        d_{\mu,ij}^{0}(\vb{x}^0) & =\langle X_{0,i}^{1}(\vb{x}^1),\alpha_{\mu}(\vb{x}^0,\vb{x}^1)X_{0,j}^{1}(\vb{x}^1)\rangle_{1}.
    \end{aligned}
\end{equation}
This yields $K_i(\Delta t,\vb{x}^0)$, and we obtain the updated low-rank factors $X_{1,i}^0(\vb{x}^0)$ by a QR decomposition $K_i(\Delta t,\vb{x}^0) = \sum_j X_{1,j}^0(\vb{x}^0) \tilde{S}_{ji}$, since the low-rank factors have to be orthogonal. This first step of the algorithm is commonly knows as the \emph{K step}. It turns out that this detour of forming $K$, instead of solving an evolution equation for $X^0(t,\vb{x}_0)$ directly, makes the algorithm robust with respect to small singular values \cite{Lubich_2014}. 

The intermediate quantity $\tilde{S}_{ij}$ is then used as an initial condition for the evolution equation of the second step of our algorithm, called the \emph{S step},
\begin{equation}\label{eq:S}
    \frac{\diff{}}{\diff{t}} S_{ij}(t)=-\sum_{k,l}S_{kl}(t)\left(e_{ijkl}-f_{ijkl}\right),
\end{equation}
with the time-independent coefficients
\begin{equation}\label{eq:S_coeff}
    \begin{aligned}
        e_{ijkl} & =\sum_\mu\langle X_{1,i}^0(\vb{x}^0)X_{0,j}^1(\vb{x}^1),\alpha_{\mu}(\vb{x}^0-\vb{\nu}^0_{\mu},\vb{x}^1-\vb{\nu}^1_{\mu})X_{1,k}^0(\vb{x}^0-\vb{\nu}^0_{\mu})X_{0,l}^1(\vb{x}^1-\vb{\nu}^1_{\mu})\rangle_{0,1},\\
        f_{ijkl} & =\sum_\mu\langle X_{1,i}^0(\vb{x}^0)X_{0,j}^1(\vb{x}^1),\alpha_{\mu}(\vb{x}^0,\vb{x}^1)X_{1,k}^{0}(\vb{x}^0)X_{0,l}^{1}(\vb{x}^1)\rangle_{0,1}.
    \end{aligned}
\end{equation}
Integrating Equation~(\ref{eq:S}) yields $\hat{S}_{ij}=S_{ij}(\Delta t)$. Note that the computation of the coefficients in the form of Equation~(\ref{eq:S_coeff}) is costly, since we have to perform inner products with respect to both partition 0 and partition 1. We can reuse the coefficients $c^0_{\mu,ij}(\vb{x}^0)$ and $d^0_{\mu,ij}(\vb{x}^0)$ for eliminating one inner product of Equation~(\ref{eq:S_coeff}), but the computational cost of the latter coefficients also scales with $\bigO(M r^2 n^0 n^1)$, where $M$ was the total number of reactions. In Section~\ref{sec:implementation} we will explain that the computation of the coefficients can be done in a more efficient way when we assume that the propensity functions are decomposable into parts that only depend on $x^0$ or on $x^1$.

In the last step of our algorithm, the \emph{L step}, we reuse the result from the S step to form $L_i(0,\vb{x}^1)=\sum_j\hat{S}_{ij} X_{0,j}^2(\vb{x}^1)$, which is then the initial condition for the evolution equation
\begin{equation}\label{eq:L}
    \partial_{t}L_{i}(t,\vb{x}^1)=\sum_\mu\sum_j\left(c_{\mu,ij}^1(\vb{x}^1) L_j(t,\vb{x}^1-\vb{\nu}^1_\mu)-d_{\mu,ij}^1(\vb{x}^1)L_j(t,\vb{x}^1)\right)
\end{equation}
with the time-independent coefficients 
\begin{equation}\label{eq:L_coeff}
    \begin{aligned}
        c_{\mu,ij}^1(\vb{x}^1) & = \langle X_{1,i}^0(\vb{x}^0),\alpha_\mu(\vb{x}^0-\vb{\nu}^0_\mu,\vb{x}^1-\vb{\nu}^1_\mu)X_{1,j}^0(\vb{x}^0-\vb{\nu}^0_\mu)\rangle_0,\\
        d_{\mu,ij}^1(\vb{x}^1) & = \langle X_{1,i}^0(\vb{x}^0),\alpha_{\mu}(\vb{x}^0,\vb{x}^1)X_{1,j}^0(\vb{x}^0)\rangle_0.
    \end{aligned}
\end{equation}
This yields $L_i(\Delta t,\vb{x}^1)$, and similar to the K step, we obtain the updated low-rank factors $X_{1,i}^1(\vb{x}^1)$  and the updated coefficient matrix $S_{1,ij}$ by a QR decomposition $L_i(\Delta t,\vb{x}^1) = \sum_j S_{1,ij}X_{1,j}^1(\vb{x}^1)$.

In total, we get an approximation for the full probability distribution at time $\Delta t$ via $P(\Delta t,\vb{x})=\sum_{i,j} X^0_{1,i}(\vb{x}^0) S_{1,ij} X^1_{1,j}(\vb{x}^1)$. An overview of the first-order projector-splitting integrator for the kinetic CME is shown in Algorithm~\ref{alg:first-order}.
\begin{algorithm}[H]
    \caption{First-order projector-splitting integrator for the kinetic CME.\label{alg:first-order}}
    \textbf{Input:} $X_{0,i}^0(\vb{x}^0),$ $S_{0,ij}$, $X_{0,j}^1(\vb{x}^1)$\\
    \textbf{Output:} $X_{1,i}^0(\vb{x}^0)$, $S_{1,ij}$, $X_{1,j}^1(\vb{x}^1)$, where $P(\Delta t, \vb{x}) \approx \sum_{i,j} X_{1,i}^0(\vb{x}^0) S_{1,ij} X_{1,j}^1(\vb{x}^1)$
    \begin{algorithmic}[1]
        \State Calculate $c^0_{\mu,ij}(\vb{x}^0)$ and $d^0_{\mu,ij}(\vb{x}^0)$ with $X_{0,i}^1(\vb{x}^1)$ using Equation~(\ref{eq:K_coeff})
        \State Integrate $K_i$ from $0$ to $\Delta t$ with initial value $K_i(0,\vb{x}^0) = \sum_j X_{0,j}^0(\vb{x}^0) S_{0,ji}$ using Equation~(\ref{eq:K})
        \State Decompose $K_i(\Delta t,\vb{x}^0) = \sum_j X_{1,j}^0(\vb{x}^0) \tilde{S}_{ji}$ via a QR factorization
        \State Calculate $e_{ijkl}$ and $f_{ijkl}$ with $X_{1,i}^0(\vb{x}^0)$, $X_{0,i}^1(\vb{x}^1)$ using Equation~(\ref{eq:S_coeff})
        \State Integrate $S_{ij}$ from $0$ to $\Delta t$ with initial value $S_{ij}(0) = \tilde{S}_{ij}$ using Equation~(\ref{eq:S}) and set $\hat{S}_{ij} = S_{ij}(\Delta t)$
        \State Calculate $c^1_{\mu,ij}(\vb{x}^1)$ and $d^1_{\mu,ij}(\vb{x}^1)$ with $X_{1,i}^0(\vb{x}^0)$ using Equation~(\ref{eq:L_coeff})
        \State Integrate $L_i$ from $0$ to $\Delta t$ with initial value $L_i(0,\vb{x}^1) = \sum_j \hat{S}_{ij} X_{0,j}^1(\vb{x}^1)$ using Equation~(\ref{eq:L})
        \State Decompose $L_i(\Delta t,\vb{x}^1) = \sum_j S_{1,ij} X_{1,j}^1(\vb{x}^1)$ via a QR factorization
    \end{algorithmic}
\end{algorithm}

\section{Tree tensor network integrator}\label{sec:ttn}
In this section we extend our matrix integrator for the kinetic CME to an integrator for binary tree tensor networks.
The main idea is to perform the partitioning recursively. Let us for example assume that the low-rank factor $X^0_i(t,\vb{x}^0)$ can be decomposed into orthonormal, time-dependent low-rank factors $X^{00}_{i_{00}}(t,\vb{x}^{00})$ and $X^{01}_{i_{01}}(t,\vb{x}^{01})$ and a connection tensor $Q^0_{ii_{00}i_{01}}(t)$ with full multilinear rank $(\tilde{r},r^{00},r^{01})$:
\begin{equation}\label{eq:DLR-X-example}
    X^0_i(t,\vb{x}^0) = \sum_{i_{00}=1}^{r^{00}} \sum_{i_{01}=1}^{r^{01}} Q^0_{ii_{00}i_{01}}(t) X^{00}_{i_{00}}(t,\vb{x}^{00}) X^{01}_{i_{01}}(t,\vb{x}^{01}),
\end{equation}
with $\vb{x}^0=(\vb{x}^{00},\vb{x}^{01})$, where $\vb{x}^{00} \in \Omega^{00}_{\vb{\zeta},\vb{\eta}} \subset \mathbb{N}^{d^{00}}_0$, $\vb{x}^{01} \in \Omega^{01}_{\vb{\zeta},\vb{\eta}} \subset \mathbb{N}^{d^{01}}_0$ and $d^{00} + d^{01} = d^0$. Similar to the matrix case, we define $n^{00}$ and $n^{01}$ to be the the size of the truncated state spaces $\Omega^{00}_{\vb{\zeta},\vb{\eta}}$ and $\Omega^{01}_{\vb{\zeta},\vb{\eta}}$, respectively. 
We observe that Equation~(\ref{eq:DLR-X-example}) is, in contrast to Equation (\ref{eq:DLR-approximation}), not a matrix decomposition, but a so-called \emph{Tucker decomposition}, since $X^0_i(t,\vb{x}^0) \equiv X^0_i(t,\vb{x}^{00},\vb{x}^{01})$ is treated as a time-dependent $r\times n^{00} \times n^{01}$~tensor.
Remember that for the matrix case the number of low-rank factors $X^0_i(t,\vb{x}^0)$ and $X^1_j(t,\vb{x}^1)$ was the same, namely $\tilde{r}$. For the Tucker decomposition the situation is different: Since the connection tensor $Q^0_{ii_{00}i_{01}}(t)$ with full multilinear rank has to obey the condition that
\begin{equation}\label{eq:Tucker-condition-example}
    (r \le r^{00} r^{01}) \land (r^{00} \le r^{01} r) \land (r^{01} \le r^{00} r),
\end{equation}
the number of low-rank factors $X^{00}_{i_{00}}(t,\vb{x}^{00})$ and $X^{01}_{i_{01}}(t,\vb{x}^{01})$ can be different from each other \cite{Ceruti_2020}.

\begin{figure}[H]
    \centering
    \begin{subfigure}[b]{0.35\textwidth} 
        \centering
        \begin{tikzpicture}
    \SetVertexStyle[TextFont=\normalsize]
    \InternalNode[label={$Q_{0 i_0 i_1}$},position=left]{root}
    \node[above=0.5cm of root](D){};
    \InternalNode[x=-1,y=-1,label={$Q^0_{i_0 i_{00} i_{01}}$},position=left]{0}
    \ExternalNode[x=1,y=-1,label={$X^1_{i_1}$},position=below]{1}
    \ExternalNode[x=-1.75,y=-2,label={$X^{00}_{i_{00}}$},position=below]{00}
    \ExternalNode[x=-0.25,y=-2,label={$X^{01}_{i_{01}}$},position=below]{01}
    
    \Edge(D)(root)
    \Edge(root)(0)
    \Edge(root)(1)
    \Edge(0)(00)
    \Edge(0)(01)
\end{tikzpicture}
        \caption{\label{fig:example-tree}}
    \end{subfigure}
    \quad
    \begin{subfigure}[b]{0.35\textwidth}
        \centering
        \begin{tikzpicture}       
    \SetVertexStyle[TextFont=\normalsize]
    \InternalNode[label={$Q^\tau_{i_\tau i_{\tau_0} i_{\tau_1}}$},position=left]{root}
    \node[above=0.5cm of root](D){};
    \ExternalNode[x=-1,y=-1,label={$X^{\tau_0}_{i_{\tau_0}}$},position=below]{0}
    \ExternalNode[x=1,y=-1,label={$X^{\tau_1}_{i_{\tau_1}}$},position=below]{1}
    
    \Edge(D)(root)
    \Edge(root)(0)
    \Edge(root)(1)
\end{tikzpicture}
        \caption{\label{fig:subtree}}
    \end{subfigure}
    \caption{(a)~Graphical representation of the right-hand side of Equation~(\ref{eq:example-tree}) as a binary tree. (b)~Graphical representation of an internal node $\tau$ (gray) with child nodes $\tau_0$ and $\tau_1$ (blue). The child nodes can either be leaves, where the actual low-rank factors $X^{\tau_0}_{i_{\tau_0}}$ and $X^{\tau_1}_{i_{\tau_1}}$ are stored, or they are also internal nodes, and the low-rank factors are recursively defined by Equation~(\ref{eq:DLR-X}).}
\end{figure}
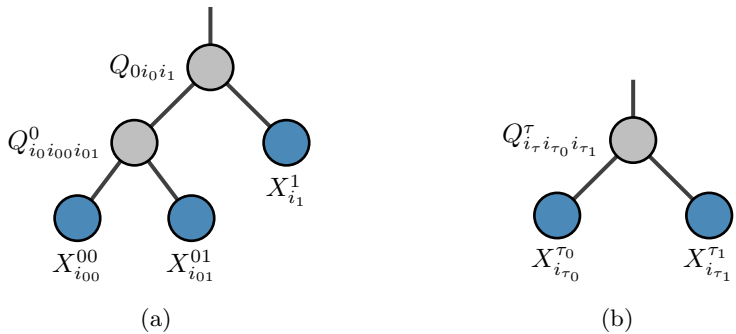

By means of Equations~(\ref{eq:DLR-approximation}) and (\ref{eq:DLR-X-example}) we can approximate the probability distribution by a product of low-rank factors with a connection tensor $Q^0_{ii_{00}i_{01}}(t)$ and a coefficient matrix $S_{ij}$. It is now convenient to replace the coefficient matrix $S_{ij}(t)$ by a $r\times r^0 \times r^1$~connection tensor $Q_{ii_{0}i_{01}}(t)$, where $r=1$. We see from Equation~(\ref{eq:Tucker-condition-example}) that $r^{0}=r^{01} \equiv \tilde{r}$ and by replacing the indices $i$ and $j$ with $i_0$ and $i_1$, respectively, we can identify $S_{ij}(t) \equiv \sum_{i=0}^{r-1} Q_{ii_{0}i_{01}}(t)$. In the following, we will again use the short-hand notation $\sum_{i_\tau}\equiv\sum_{i_\tau=0}^{r^\tau-1}$ for a partition $\tau$. 
Thus we can write Equation~(\ref{eq:DLR-approximation}) in a similar fashion as Equation~(\ref{eq:DLR-X-example}): 
\begin{equation}
    P(t,\vb{x}) \approx \sum_{i} \sum_{i_0} \sum_{i_1} Q_{ii_0i_1}(t) X^0_{i_0}(t,\vb{x}^0) X^1_{i_1}(t,\vb{x}^1).
    \label{eq:DLR-approximation-new}
\end{equation}
Then the approximation of the full probability distribution via Equations~(\ref{eq:DLR-approximation-new}) and (\ref{eq:DLR-X-example}) reads as
\begin{equation}\label{eq:example-tree}
    \begin{aligned}
        P(t,\vb{x}) \approx & \sum_{i}\sum_{i_0}\sum_{i_1}\sum_{i_{00}}\sum_{i_{01}}Q_{ii_{0}i_{1}}(t)Q_{i_{0}i_{00}i_{01}}^{0}(t) X_{i_{00}}^{00}(t,\vb{x}^{00}) X_{i_{01}}^{01}(t,\vb{x}^{01}) X_{i_{1}}^{1}(t,\vb{x}^1).
    \end{aligned}
\end{equation}
We can interpret the right-hand side of this equation graphically as follows: For the connection tensors, we associate incoming edges with the first index and outgoing edges with the remaining two indices, whereas for the low-rank factors we identify the index with an incoming edge. As can be seen from Figure~\ref{fig:example-tree}, we end up with a binary tree structure, where the coefficient tensors are situated at the internal nodes (gray) and the low-rank factors are at the leaves of the tree (blue). Hence, Equation~(\ref{eq:example-tree}) is an example of the so-called \emph{tree tensor network} (TTN) decomposition.

We will now generalize this to arbitrary binary TTNs. We assume that the low-rank factors can be recursively decomposed as
\begin{equation}\label{eq:DLR-X}
    X^\tau_{i_\tau}(t,\vb{x}^\tau) = \sum_{i_{\tau_0}} \sum_{i_{\tau_1}} Q^\tau_{i_\tau i_{\tau_0}i_{\tau_1}}(t) X^{\tau_0}_{i_{\tau_0}}(t,\vb{x}^{\tau_0}) X^{\tau_1}_{i_{\tau_1}}(t,\vb{x}^{\tau_1}),
\end{equation}
for a given partition $\tau$ with subpartitions $\tau_0$ and $\tau_1$, such that $\vb{x}^\tau = (\vb{x}^{\tau_0}, \vb{x}^{\tau_1})$. Therefore the approximation of the probability distribution given by Equation~(\ref{eq:DLR-approximation-new}) can also be expressed by the low-rank factor of the root node:
\begin{equation*}
    P(t,\vb{x})\approx \sum_i X_i(t, \vb{x}) = X_0(t,\vb{x}).
\end{equation*}
Moreover, we demand that the low-rank factors $X^{\tau_0}_{i_{\tau_0}}(t,\vb{x}^{\tau_0})$ and $X^{\tau_1}_{i_{\tau_1}}(t,\vb{x}^{\tau_1})$ again have to be orthonormal and the connection tensors $Q^{\tau}_{i_\tau i_{\tau_0} i_{\tau_1}}$ have full multilinear rank $(\tilde{r^\tau},r^{\tau_0},r^{\tau_1})$ and obey the general form of Equation~(\ref{eq:Tucker-condition-example}),
\begin{equation}\label{eq:Tucker-condition}
    (r^\tau \le r^{\tau_0} r^{\tau_1}) \land (r^{\tau_0} \le r^{\tau_1} r) \land (r^{\tau_1} \le r^{\tau_0} r).
\end{equation}
In context of the graphical representation of binary trees, $\tau_0$ and $\tau_1$ are subtrees of the tree $\tau$. Theoretically, we could do this decomposition until the low-rank factors on the leaves of the resulting tree only depend on a single physical dimension, i.e., a single species. It would be even possible to decompose a physical dimension into virtual dimensions similar to the quantized tensor train decomposition. In our view it is however often preferable to not fully decompose the probability distribution, but to keep species, which are coupled by important reaction pathways, together in a single partition. We emphasize that the structure of the tree is chosen based on the specific problem and during the time integration the tree structure is fixed.
The approximation of the probability distribution with Equations~(\ref{eq:DLR-approximation-new}) and (\ref{eq:DLR-X}) is now fully described by low-rank factors on the leaves and by connection tensors on the internal nodes of the resulting binary tree. We do not store the low-rank factors on the internal nodes, but compute them recursively via Equation~(\ref{eq:DLR-X}). For the time evolution of $P(t,\vb{x})$ via the CME we have to derive an integrator which advances the low-rank factors and the connection tensors in time without forming the full probability distribution.

\subsection{Algorithm}\label{sec:ttn-algorithm}
We will now describe the algorithm for the projector-splitting tree tensor network (PS-TTN) integrator for the kinetic chemical master equation (for more details on the integrator we refer the reader to \cite{Ceruti_2020}).

The replacement of the coefficient matrix by a $r \times r^0 \times r^1$ connection tensor (with $r=1$) on the root of the binary tree had the effect that an internal node $\tau$ (including the root) is associated with a connection tensor $Q^\tau_{i_\tau i_{\tau_0} i_{\tau_1}}(t)$ obeying the rank condition (\ref{eq:Tucker-condition}). The two child nodes of $\tau$ are in turn associated with low-rank factors $X^{\tau_0}_{i_{\tau_0}}(t,\vb{x}^{\tau_0})$ and $X^{\tau_1}_{i_{\tau_1}}(t,\vb{x}^{\tau_1})$ (see Figure~\ref{fig:subtree}). If the child nodes or are again internal nodes, then the low-rank factors are recursively defined via Equation~(\ref{eq:DLR-X}). Only if $\tau_0$ or $\tau_1$ are leaves, then the corresponding low-rank factors are actually stored.

We assume that the initial condition of the probability distribution is approximated by 
\begin{equation}
    P(0,\vb{x}) \approx \sum_{i} \sum_{i_0} \sum_{i_1} Q_{0,ii_0i_1} X^0_{0,i_0}(\vb{x}^0) X^1_{0,i_1}(\vb{x}^1)
\end{equation}
and for the low-rank factors recursively by
\begin{equation}\label{eq:DLR-X-initial}
    X^\tau_{0,i_\tau}(\vb{x}^\tau) = \sum_{i_{\tau_0}} \sum_{i_{\tau_1}} Q^\tau_{0,i_\tau i_{\tau_0}i_{\tau_1}} X^{\tau_0}_{0,i_{\tau_0}}(\vb{x}^{\tau_0}) X^{\tau_1}_{0,i_{\tau_1}}(\vb{x}^{\tau_1}).
\end{equation}
Similar to the matrix case, we can either have the situation where $P(t,\vb{x})$ is explicitly given in the TTN format and the connection tensors and low-rank factors can be directly inferred. Otherwise, the probability distribution has to be approximated with a (randomized) truncated higher-order SVD.

Our aim now is to advance the internal node $\tau$ in time by a single time step $\Delta t$. The general procedure is similar to the matrix integrator, where updating the left and right child nodes $\tau_0$ and $\tau_1$ of Figure~\ref{fig:subtree} corresponds to the K and L step, respectively. For each of the child nodes we also have to perform an S step. In the first two steps we advance the low-rank factors in the left and right child nodes $\tau_0$ and $\tau_1$ in time, and in a third step we update a quantity closely related to the connection tensor of the internal node $\tau$.

\subsubsection{First step}
We start with updating the low-rank factor of the left child node, $X^{\tau_0}_{0,i_{\tau_0}}(\vb{x}^{\tau_0})$. Since Equation~(\ref{eq:DLR-X}) is not a matrix, but a Tucker decomposition, we now have to apply the Tucker tensor integrator of \cite{Lubich_2018} for the time integration of $\tau_0$. 

Therefore, we first have to perform a QR decomposition of the tensor $C^\tau_{0,i_\tau i_{\tau_0} i_{\tau_1}}$, which is related to the connection tensor $Q^\tau_{0,i_\tau i_{\tau_0} i_{\tau_1}}$: For the root, we set $C_{0,i i_0 i_1} = Q_{0,i i_0 i_1}$ and we will see that for every internal node we have already computed the required tensor in the parent node. Therefore we can assume that we know this quantity, and we can perform a QR decomposition by treating this tensor as a matrix $\tilde{C}^\tau_{0,\beta_{\tau_0} i_{\tau_0}} \equiv C^\tau_{0,i_\tau i_{\tau_0} i_{\tau_1}}$ with a combined index $\beta_{\tau_0}=i_\tau + r^\tau i_{\tau_1}$: 
\begin{equation*}
    C^\tau_{0,i_\tau i_{\tau_0} i_{\tau_1}} = \sum_{j_{\tau_0}} G^{\tau_0}_{i_\tau j_{\tau_0} i_{\tau_1}} S^{\tau_0}_{0,i_{\tau_0} j_{\tau_0}}.
\end{equation*}

Next, we have to determine whether $\tau_0$ is a leaf node or an internal node.
If the latter is the case, we compute $C^{\tau_0}_{0,i_{\tau_0} i_{\tau_{00}} i_{\tau_{01}}} = \sum_{j_{\tau_0}} Q^{\tau_0}_{0,j_{\tau_0} i_{\tau_{00}} i_{\tau_{01}}} S^{\tau_0}_{0,j_{\tau_0} i_{\tau_0}}$ for the child node and recursively do the integration for the internal node $\tau_0$ which yields, as we will see in the third step, an updated $C^{\tau_0}_{3,i_{\tau_0} i_{\tau_{00}} i_{\tau_{01}}}$ and updated low-rank factors $X^{\tau_{00}}_{1,i_{\tau_{00}}}(\vb{x}^{\tau_{00}})$ and $X^{\tau_{01}}_{1,i_{\tau_{01}}}(\vb{x}^{\tau_{01}})$. After that, we perform the QR decomposition
\begin{equation*}
    C^{\tau_0}_{3,i_{\tau_0} i_{\tau_{00}} i_{\tau_{01}}} = \sum_{j_{\tau_0}} Q^{\tau_0}_{1,j_{\tau_0} i_{\tau_{00}} i_{\tau_{01}}} \tilde{S}^{\tau_0}_{j_{\tau_0} i_{\tau_0}},
\end{equation*}
which can again be achieved by identifying the tensor with a matrix $\hat{C}^{\tau_0}_{\gamma_{\tau_0} i_{\tau_0}} \equiv C^{\tau_0}_{1,i_{\tau_0} i_{\tau_{00}} i_{\tau_{01}}}$ using a combined index $\gamma_{\tau_0}=i_{\tau_{00}} + r^{\tau_{00}} i_{\tau_{01}}$. This yields the updated connection tensor $Q^{\tau_0}_{1,i_{\tau_0} i_{\tau_{00}} i_{\tau_{01}}}$ and an intermediate quantity $\tilde{S}^{\tau_0}_{i_{\tau_0} j_{\tau_0}}$, which will be used later. The updated low-rank factor at $\tau_0$ is given by 
\begin{equation}\label{eq:DLR-X-update}
    X^{\tau_0}_{1,i_{\tau_0}}(\vb{x}_{\tau_0}) = \sum_{i_{\tau_{00}}} \sum_{i_{\tau_{01}}} Q^{\tau_0}_{1,i_{\tau_0} i_{\tau_{00}} i_{\tau_{01}}} X^{\tau_{00}}_{1,i_{\tau_{00}}}(\vb{x}_{\tau_{00}}) X^{\tau_{01}}_{1,i_{\tau_{01}}}(\vb{x}_{\tau_{01}}).
\end{equation}

If $\tau_0$ is on the other hand a leaf node, we compute $W^{\tau_0}_{i_\tau i_{\tau_0}}(\vb{x}^{\tau_1}) = \sum_{i_{\tau_1}} G^{\tau_0}_{i_\tau i_{\tau_0} i_{\tau_1}} X^{\tau_1}_{0,i_{\tau_1}}(\vb{x}^{\tau_1})$. Note that $X^{\tau_1}_{0,i_{\tau_1}}(\vb{x}^{\tau_1})$ is directly stored when $\tau_1$ is a leaf node, or it can be calculated by Equation~(\ref{eq:DLR-X-initial}) if $\tau_1$ is an internal node. Then, we use $K^{\tau_0}_{i_{\tau_0}}(0,\vb{x}^{\tau_0})=\sum_{j_{\tau_0}} X^{\tau_0}_{0,j_{\tau_0}}(\vb{x}^{\tau_0}) S^{\tau_0}_{0,j_{\tau_0} i_{\tau_0}}$ as an initial value for the integration of
\begin{equation}\label{eq:K0-TTN}
    \partial_{t}K^{\tau_0}_{i_{\tau_0}}(t,\vb{x}^{\tau_0})=\sum_\mu\sum_{j_{\tau_0}}\left(c_{\mu,i_{\tau_0} j_{\tau_0}}^{\tau_0}(\vb{x}^{\tau_0}) K_{j_{\tau_0}}(t,\vb{x}^{\tau_0}-\vb{\nu}^{\tau_0}_{\mu})-d_{\mu,i_{\tau_0}j_{\tau_0}}^{\tau_0}(\vb{x}^{\tau_0})K_{j_{\tau_0}}(t,\vb{x}^{\tau_0})\right)
\end{equation}
from $0$ to $\Delta t$. The time-independent coefficients are recursively defined by
\begin{equation}\label{eq:K0-TTN_coeff}
    \begin{aligned}
        c_{\mu,i_{\tau_0}j_{\tau_0}}^{\tau_0}(\vb{x}^{\tau_0}) & =\sum_{i_\tau,j_\tau}\langle W_{i_\tau i_{\tau_0}}^{\tau_0}(\vb{x}^{\tau_1}), c_{\mu,i_{\tau}j_{\tau}}^{\tau}(\vb{x}^{\tau}) W_{j_\tau j_{\tau_0}}^{\tau_0}(\vb{x}^{\tau_1}-\vb{\nu}^{\tau_1}_{\mu})\rangle_{\tau_1},\\
        d_{\mu,i_{\tau_0}j_{\tau_0}}^{\tau_0}(\vb{x}^{\tau_0}) & =\sum_{i_\tau,j_\tau}\langle W_{i_\tau i_{\tau_0}}^{\tau_0}(\vb{x}^{\tau_1}), d_{\mu,i_{\tau}j_{\tau}}^{\tau}(\vb{x}^{\tau}) W_{j_\tau j_{\tau_0}}^{\tau_0}(\vb{x}^{\tau_1})\rangle_{\tau_1},
    \end{aligned}
\end{equation}
with the initial condition
\begin{equation*}
    \begin{aligned}
        c_{\mu,ij}^{0}(\vb{x}^0) & =\langle X_{0,i}^{1}(\vb{x}^1),\alpha_{\mu}(\vb{x}^0-\vb{\nu}^0_{\mu},\vb{x}^1-\vb{\nu}^1_{\mu})X_{0,j}^1(\vb{x}^1-\vb{\nu}^1_{\mu})\rangle_{1},\\
        d_{\mu,ij}^{0}(\vb{x}^0) & =\langle X_{0,i}^{1}(\vb{x}^1),\alpha_{\mu}(\vb{x}^0,\vb{x}^1)X_{0,j}^{1}(\vb{x}^1)\rangle_{1}.
    \end{aligned}
\end{equation*}
We then perform the QR decomposition
\begin{equation}\label{eq:K0-TTN-QR}
    K^{\tau_0}_{i_{\tau_0}}(\Delta t, \vb{x}^{\tau_0}) = \sum_{j_{\tau_0}} X^{\tau_0}_{1,j_{\tau_0}}(\vb{x}^{\tau_0}) \tilde{S}^{\tau_0}_{j_{\tau_0} i_{\tau_0}},
\end{equation}
yielding the updated low-rank factor on the leaf node $\tau_0$ and an intermediate quantity $\tilde{S}^{\tau_0}_{i_{\tau_0} j_{\tau_0}}$.

In both situations, $\tau_0$ being a leaf node or an internal node, we now use the intermediate $\tilde{S}^{\tau_0}_{i_{\tau_0} j_{\tau_0}}$ as an initial value for the evolution equation
\begin{equation}\label{eq:S0-TTN}
    \frac{\diff{}}{\diff{t}}S^{\tau_0}_{i_{\tau_0} j_{\tau_0}}(t) = -\sum_{k_{\tau_0},l_{\tau_0}} S^{\tau_0}_{k_{\tau_0} l_{\tau_0}}(t) \left( e^{\tau_0}_{i_{\tau_0} j_{\tau_0} k_{\tau_0} l_{\tau_0}} - f^{\tau_0}_{i_{\tau_0} j_{\tau_0} k_{\tau_0} l_{\tau_0}} \right),
\end{equation}
with the time-independent coefficients
\begin{equation}\label{eq:S0-TTN_coeff}
    \begin{aligned}
        e^{\tau_0}_{i_{\tau_0} j_{\tau_0} k_{\tau_0} l_{\tau_0}} & =\sum_{\mu} \langle X^{\tau_0}_{1,i_{\tau_0}}(\vb{x}^{\tau_0}), c_{\mu,j_{\tau_0}l_{\tau_0}}^{\tau_0}(\vb{x}^{\tau_0}), X_{1,k_{\tau_0}}^{\tau_0}(\vb{x}^{\tau_0}-\vb{\nu}_{\mu}^{\tau_0})\rangle_{\tau_0}\\
        f^{\tau_0}_{i_{\tau_0} j_{\tau_0} k_{\tau_0} l_{\tau_0}} & =\sum_{\mu} \langle X^{\tau_0}_{1,i_{\tau_0}}(\vb{x}^{\tau_0}), d_{\mu,j_{\tau_0}l_{\tau_0}}^{\tau_0}(\vb{x}^{\tau_0}), X_{1,k_{\tau_0}}^{\tau_0}(\vb{x}^{\tau_0})\rangle_{\tau_0},
    \end{aligned}
\end{equation}
which is similar to the S step of the matrix integrator. We integrate this equation from $0$ to $\Delta t$ and set $S^{\tau_0}_{1,i_{\tau_0} j_{\tau_0}} = S^{\tau_0}_{i_{\tau_0} j_{\tau_0}}(\Delta t)$. Finally, we compute $C^{\tau}_{1,i_\tau i_{\tau_0} i_{\tau_1}} = \sum_{j_{\tau_0}} G^{\tau_0}_{i_\tau j_{\tau_0} i_{\tau_1}} S^{\tau_0}_{1,i_{\tau_0} j_{\tau_0}}$, which is the starting point for updating the right child node $\tau_1$.

Summing up, we have updated in this step the low-rank factor at the left child node to $X^{\tau_0}_{1,i_{\tau_0}}(\vb{x}_{\tau_0})$. When $\tau_0$ is an internal node, it can be computed by Equation~(\ref{eq:DLR-X-update}), where the quantities on the right-hand side are determined by the recursive application of the integrator.

\subsubsection{Second step}
In the second step, we update the low-rank factor of the right child, $X^{\tau_1}_{0,i_{\tau_1}}(\vb{x}^{\tau_1})$. This is very similar to the update of the left child, but here we start with the QR decomposition
\begin{equation*}
    C^\tau_{1,i_\tau i_{\tau_0} i_{\tau_1}} = \sum_{j_{\tau_1}} G^{\tau_1}_{i_\tau i_{\tau_0} j_{\tau_1}} S^{\tau_1}_{0,i_{\tau_1} j_{\tau_1}},
\end{equation*}
which can be achieved by writing the tensor as a matrix $\tilde{C}^\tau_{1,\beta_{\tau_1} i_{\tau_1}} \equiv C^\tau_{1,i_\tau i_{\tau_0} i_{\tau_1}}$, with a combined index $\beta_{\tau_1}=i_\tau + r^\tau i_{\tau_0}$.

Now we again have to determine whether $\tau_0$ is a leaf node or an internal node.
If it is an internal node, we compute $C^{\tau_1}_{0,i_{\tau_1} i_{\tau_{10}} i_{\tau_{11}}} = \sum_{j_{\tau_1}} Q^{\tau_1}_{0,j_{\tau_1} i_{\tau_{10}} i_{\tau_{11}}} S^{\tau_1}_{0,j_{\tau_1} i_{\tau_1}}$ for the child node and recursively do the integration for the internal node $\tau_1$, which yields an updated $C^{\tau_1}_{3,i_{\tau_1} i_{\tau_{10}} i_{\tau_{11}}}$ and updated low-rank factors $X^{\tau_{10}}_{1,i_{\tau_{10}}}(\vb{x}^{\tau_{10}})$ and $X^{\tau_{11}}_{1,i_{\tau_{11}}}(\vb{x}^{\tau_{11}})$. Next, we perform the QR decomposition
\begin{equation*}
    C^{\tau_1}_{3,i_{\tau_1} i_{\tau_{10}} i_{\tau_{11}}} = \sum_{j_{\tau_1}} Q^{\tau_1}_{1,j_{\tau_1} i_{\tau_{10}} i_{\tau_{11}}} \tilde{S}^{\tau_1}_{j_{\tau_1} i_{\tau_1}},
\end{equation*}
by identifying the tensor with a matrix $\hat{C}^{\tau_1}_{\gamma_{\tau_1} i_{\tau_1}} \equiv C^{\tau_1}_{3,i_{\tau_1} i_{\tau_{10}} i_{\tau_{11}}}$ using a combined index $\gamma_{\tau_1}=i_{\tau_{10}} + r^{\tau_{10}} i_{\tau_{11}}$. This yields an update for the connection tensor $Q^{\tau_1}_{1,i_{\tau_1} i_{\tau_{10}} i_{\tau_{11}}}$ and the intermediate quantity $\tilde{S}^{\tau_1}_{i_{\tau_1} j_{\tau_1}}$.

If $\tau_1$ is a leaf node, we compute $W^{\tau_1}_{i_\tau i_{\tau_1}}(\vb{x}^{\tau_0}) = \sum_{i_{\tau_0}} G^{\tau_1}_{i_\tau i_{\tau_0} i_{\tau_1}} X^{\tau_0}_{1,i_{\tau_0}}(\vb{x}^{\tau_0})$. Recall that $X^{\tau_0}_{1,i_{\tau_0}}(\vb{x}^{\tau_0})$ is directly stored when $\tau_0$ is a leaf node, or it can be computed via Equation~(\ref{eq:DLR-X-update}) if $\tau_0$ is an internal node. Then, we use $K^{\tau_1}_{i_{\tau_1}}(0,\vb{x}^{\tau_1})=\sum_{j_{\tau_1}} X^{\tau_1}_{0,j_{\tau_1}}(\vb{x}^{\tau_1}) S^{\tau_1}_{0,j_{\tau_1} i_{\tau_1}}$ as an initial value for the integration of
\begin{equation}\label{eq:K1-TTN}
    \partial_{t}K^{\tau_1}_{i_{\tau_1}}(t,\vb{x}^{\tau_1})=\sum_\mu\sum_{j_{\tau_1}}\left(c_{\mu,i_{\tau_1} j_{\tau_1}}^{\tau_1}(\vb{x}^{\tau_1}) K_{j_{\tau_1}}(t,\vb{x}^{\tau_1}-\vb{\nu}^{\tau_1}_{\mu})-d_{\mu,i_{\tau_1}j_{\tau_1}}^{\tau_1}(\vb{x}^{\tau_1})K_{j_{\tau_1}}(t,\vb{x}^{\tau_1})\right)
\end{equation}
from $0$ to $\Delta t$. The time-independent coefficients $c_{\mu,i_{\tau_1} j_{\tau_1}}^{\tau_1}(\vb{x}^{\tau_1})$ and $d_{\mu,i_{\tau_1} j_{\tau_1}}^{\tau_1}(\vb{x}^{\tau_1})$ are similar to the coefficients for the leaf node $\tau_0$ defined in Equation~(\ref{eq:K0-TTN_coeff}).
We perform the QR decomposition
\begin{equation*}
    K^{\tau_1}_{i_{\tau_1}}(\Delta t, \vb{x}^{\tau_1}) = \sum_{j_{\tau_1}} X^{\tau_1}_{1,j_{\tau_1}}(\vb{x}^{\tau_1}) \tilde{S}^{\tau_1}_{j_{\tau_1} i_{\tau_1}},
\end{equation*}
which yields the updated low-rank factor on the leaf node $\tau_1$ and an intermediate quantity $\tilde{S}^{\tau_1}_{i_{\tau_1} j_{\tau_1}}$.

As for the left child node $\tau_0$, we use the intermediate $\tilde{S}^{\tau_0}_{i_{\tau_0} j_{\tau_0}}$ in both cases ($\tau_0$ being a leaf node or an internal node) as an initial value for the evolution equation
\begin{equation}\label{eq:S1-TTN}
    \frac{\diff{}}{\diff{t}}S^{\tau_1}_{i_{\tau_1} j_{\tau_1}}(t) = -\sum_{k_{\tau_1},l_{\tau_1}} S^{\tau_1}_{k_{\tau_1} l_{\tau_1}}(t) \left( e^{\tau_1}_{i_{\tau_1} j_{\tau_1} k_{\tau_1} l_{\tau_1}} - f^{\tau_1}_{i_{\tau_1} j_{\tau_1} k_{\tau_1} l_{\tau_1}} \right),
\end{equation}
where the time-independent coefficients $e^{\tau_1}_{i_{\tau_1} j_{\tau_1} k_{\tau_1} l_{\tau_1}}$ and $f^{\tau_1}_{i_{\tau_1} j_{\tau_1} k_{\tau_1} l_{\tau_1}}$ are similar to the coefficients for the left child node $\tau_0$ defined in Equation~(\ref{eq:S0-TTN_coeff}). We integrate this equation from $0$ to $\Delta t$ and set $S^{\tau_1}_{1,i_{\tau_1} j_{\tau_1}} = S^{\tau_1}_{i_{\tau_1} j_{\tau_1}}(\Delta t)$. Finally, we compute $C^{\tau}_{2,i_\tau i_{\tau_0} i_{\tau_1}} = \sum_{j_{\tau_1}} G^{\tau_1}_{i_\tau i_{\tau_0} j_{\tau_1}} S^{\tau_1}_{1,i_{\tau_1} j_{\tau_1}}$, which is used as an initial condition for the evolution equation in the third step.

In summary, we have updated in this step the low-rank factors for the right child node $\tau_1$ to $X^{\tau_1}_{1,i_{\tau_1}}(\vb{x}_{\tau_1})$. If $\tau_1$ is an internal node, we can compute them by
\begin{equation}
    X^{\tau_1}_{1,i_{\tau_1}}(\vb{x}_{\tau_1}) = \sum_{i_{\tau_{10}}} \sum_{i_{\tau_{11}}} Q^{\tau_1}_{1,i_{\tau_1} i_{\tau_{10}} i_{\tau_{11}}} X^{\tau_{10}}_{1,i_{\tau_{10}}}(\vb{x}_{\tau_{10}}) X^{\tau_{11}}_{1,i_{\tau_{11}}}(\vb{x}_{\tau_{11}}),
\end{equation}
where the quantities on the right-hand side are determined by applying the integrator recursively.

\subsubsection{Third step}
In the third and last step, we use $C^{\tau}_{2,i_\tau i_{\tau_0} i_{\tau_1}}$ as an initial condition for the time integration of
\begin{equation}\label{eq:C-TTN}
    \frac{\diff{}}{\diff{t}}C^{\tau}_{i_\tau i_{\tau_0} i_{\tau_1}}(t) = \sum_{j_\tau} \sum_{j_{\tau_0}} \sum_{j_{\tau_1}} C^{\tau}_{j_\tau j_{\tau_0} j_{\tau_1}}(t) \left( g^\tau_{i_\tau i_{\tau_0} i_{\tau_1} j_\tau j_{\tau_0} j_{\tau_1}} - h^\tau_{i_\tau i_{\tau_0} i_{\tau_1} j_\tau j_{\tau_0} j_{\tau_1}} \right)
\end{equation}
from $0$ to $\Delta t$, with the time-independent coefficients
\begin{equation}\label{eq:C-TTN_coeff}
    \begin{aligned}
        g^{\tau}_{i_\tau i_{\tau_0} i_{\tau_1} j_\tau j_{\tau_0} j_{\tau_1}} & = \sum_{\mu} \langle X^{\tau_0}_{1,i_{\tau_0}}(\vb{x}^{\tau_0}) X^{\tau_1}_{1,i_{\tau_1}}(\vb{x}^{\tau_1}), c_{\mu,i_\tau j_\tau}^\tau(\vb{x}^\tau) X_{1,j_{\tau_0}}^{\tau_0}(\vb{x}^{\tau_0}-\vb{\nu}_{\mu}^{\tau_0}) X_{1,j_{\tau_1}}^{\tau_1}(\vb{x}^{\tau_1}-\vb{\nu}_{\mu}^{\tau_1})\rangle_{\tau_0,\tau_1}\\
        h^{\tau}_{i_\tau i_{\tau_0} i_{\tau_0} j_\tau j_{\tau_1} j_{\tau_1}} &= \sum_{\mu} \langle X^{\tau_0}_{1,i_{\tau_0}}(\vb{x}^{\tau_0}) X^{\tau_1}_{1,i_{\tau_1}}(\vb{x}^{\tau_1}), d_{\mu,i_\tau j_\tau}^\tau(\vb{x}^\tau) X_{1,j_{\tau_0}}^{\tau_0}(\vb{x}^{\tau_0}) X_{1,j_{\tau_1}}^{\tau_1}(\vb{x}^{\tau_1})\rangle_{\tau_0,\tau_1}.
    \end{aligned}
\end{equation}
We finally set $C^{\tau}_{3, i_\tau i_{\tau_0} i_{\tau_1}}=C^{\tau}_{i_\tau i_{\tau_0} i_{\tau_1}}(\Delta t)$ or, if $\tau$ is the root, $Q_{1,i i_0 i_1} = C_{i i_0 i_1}(\Delta t)$. Note that there is a relation between the coefficients for Equation~(\ref{eq:S1-TTN}) and the coefficients~(\ref{eq:C-TTN_coeff}), namely
\begin{equation}\label{eq:ef-gh_coeff}
    \begin{aligned}
        e^{\tau_1}_{i_{\tau_1} k_{\tau_1} j_{\tau_1} l_{\tau_1}} &= \sum_{i_\tau, j_\tau} \sum_{i_{\tau_0}, j_{\tau_0}} G^{\tau_1}_{i_\tau i_{\tau_0} k_{\tau_1}} G^{\tau_1}_{j_\tau j_{\tau_0} l_{\tau_1}} g^{\tau}_{i_\tau i_{\tau_0} i_{\tau_1} j_\tau j_{\tau_0} j_{\tau_1}},\\
        f^{\tau_1}_{i_{\tau_1} k_{\tau_1} j_{\tau_1} l_{\tau_1}} &= \sum_{i_\tau, j_\tau} \sum_{i_{\tau_0}, j_{\tau_0}} G^{\tau_1}_{i_\tau i_{\tau_0} k_{\tau_1}} G^{\tau_1}_{j_\tau j_{\tau_0} l_{\tau_1}} h^{\tau}_{i_\tau i_{\tau_0} i_{\tau_1} j_\tau j_{\tau_0} j_{\tau_1}}.
    \end{aligned}
\end{equation}
We use these two identities in the unit tests of our implementation for verification purposes.
Combining the three steps and starting at the root node, the integrator climbs up and down the tree and updates thereby the low-rank factors and coefficient tensors. Remember that for internal nodes, the low-rank factors are not stored directly, but they have to be computed on the fly for the evaluation of the coefficients (\ref{eq:K0-TTN_coeff}), (\ref{eq:S0-TTN_coeff}) and (\ref{eq:C-TTN_coeff}). Moreover, the coefficients (\ref{eq:K0-TTN_coeff}) are recursively defined, which also complicates the calculation. We will address the efficient computation of the coefficients in more detail in the next section.

\subsection{Implementation}\label{sec:implementation}
In order to compute the coefficients for the evolution equations efficiently, we demand that the propensity functions $\alpha_\mu(\vb{x})$ satisfy the \emph{factorization property} at the root node,
\begin{equation}\label{eq:factorization}
    \alpha_\mu(\vb{x}) = \alpha^0_\mu(\vb{x}^0) \alpha^1_\mu(\vb{x}^1),
\end{equation}
and at all other internal nodes $\tau$:
\begin{equation}\label{}
    \alpha^{\tau}_\mu(\vb{x}^\tau) = \alpha^{\tau_0}_\mu(\vb{x}^{\tau_0}) \alpha^{\tau_1}_\mu(\vb{x}^{\tau_1}).
\end{equation}
The factorization property is ubiquitous in most biological systems: It is fulfilled by elementary reactions (see, for example, \cite{Erban_2020}) and it is also valid for propensities stemming from effective models such as Michaelis-Menten kinetics.

Let us first consider the coefficients for the first step of Section~\ref{sec:ttn-algorithm}, where the left child node $\tau_0$ is updated. For computing the coefficients, we have to calculate inner products between propensity functions and low-rank factors. We therefore set
\begin{equation}\label{eq:A_coeff}
    \begin{aligned}
        A^\tau_{\mu,i_\tau j_\tau} &= \langle X^\tau_{0,i_\tau}(\vb{x}^\tau), \alpha^\tau_\mu(\vb{x}^\tau-\vb{\nu}^\tau) X^\tau_{0,j_\tau}(\vb{x}^\tau-\vb{\nu}^\tau) \rangle_\tau, \\
        B^\tau_{\mu,i_\tau j_\tau} &= \langle X^\tau_{0,i_\tau}(\vb{x}^\tau), \alpha^\tau_\mu(\vb{x}^\tau) X^\tau_{0,j_\tau}(\vb{x}^\tau) \rangle_\tau
    \end{aligned}
\end{equation}
with the initial conditions $X^{\tau}_{0,i_\tau}(\vb{x}^\tau)$ for the low-rank factors. Before starting with the time evolution of the TTN, we compute and store $A^\tau_{\mu,i_\tau j_\tau}$ and $B^\tau_{\mu,i_\tau j_\tau}$ for all nodes. For the leaves we can compute them directly and for the internal nodes this can be done efficiently by exploiting the factorization property, which yields
\begin{equation}\label{eq:A_coeff-recursion}
    \begin{aligned}
        A^\tau_{\mu,i_{\tau} j_{\tau}} &= \sum_{i_{\tau_0},j_{\tau_0}} \sum_{i_{\tau_1},j_{\tau_1}} Q^\tau_{i_\tau i_{\tau_0} i_{\tau_1}} Q^\tau_{j_\tau j_{\tau_0} j_{\tau_1}} A^{\tau_0}_{\mu,i_{\tau_0} j_{\tau_0}} A^{\tau_1}_{\mu,i_{\tau_1} j_{\tau_1}},\\
        B^\tau_{\mu,i_{\tau} j_{\tau}} &= \sum_{i_{\tau_0},j_{\tau_0}} \sum_{i_{\tau_1},j_{\tau_1}} Q^\tau_{i_\tau i_{\tau_0} i_{\tau_1}} Q^\tau_{j_\tau j_{\tau_0} j_{\tau_1}} B^{\tau_0}_{\mu,i_{\tau_0} j_{\tau_0}} B^{\tau_1}_{\mu,i_{\tau_1} j_{\tau_1}}.
    \end{aligned}
\end{equation}
We emphasize that we have to do this only once before the simulation and we update a coefficient with Equation~(\ref{eq:A_coeff}) or (\ref{eq:A_coeff-recursion}) only when the low-rank factor or the connection tensor of the corresponding node has been updated.

The complexity for calculating a single coefficient at a leaf node, where the low-rank factor is stored, is $\bigO(Mn^\tau(r^\tau)^2)$, whereas at an internal node, where we use the recursion relation (\ref{eq:A_coeff-recursion}), the computational cost scales with $\bigO(M (\bar{r}^\tau)^4)$, where $\bar{r}^\tau = \max{(r^\tau,r^{\tau_0},r^{\tau_1})}$ and $M$ is the number of reactions.

Due to the factorization property, we can now replace $c^{\tau_0}_{\mu,i_{\tau_0} j_{\tau_0}}(\vb{x}^{\tau_0})$ and $d^{\tau_0}_{\mu,i_{\tau_0} j_{\tau_0}}(\vb{x}^{\tau_0})$ of Equation~(\ref{eq:K0-TTN_coeff}) for the K step with coefficients which show no population number dependency:
\begin{equation}\label{eq:K0-TTN_coeff-new}
    \begin{aligned}
        a^{\tau_0}_{\mu, i_{\tau_0} j_{\tau_0}} &= \sum_\mu \sum_{i_\tau, j_\tau} \sum_{i_{\tau_1}, j_{\tau_1}} G^{\tau_0}_{i_\tau i_{\tau_0} i_{\tau_1}} G^{\tau_0}_{j_\tau j_{\tau_0} j_{\tau_1}} A^{\tau_1}_{\mu,i_{\tau_1}j_{\tau_1}} a^\tau_{\mu, i_\tau j_\tau}, \\
        b^{\tau_0}_{\mu, i_{\tau_0} j_{\tau_0}} &= \sum_\mu \sum_{i_\tau, j_\tau} \sum_{i_{\tau_1}, j_{\tau_1}} G^{\tau_0}_{i_\tau i_{\tau_0} i_{\tau_1}} G^{\tau_0}_{j_\tau j_{\tau_0} j_{\tau_1}} B^{\tau_1}_{\mu,i_{\tau_1} j_{\tau_1}} b^\tau_{\mu, i_\tau j_\tau},
    \end{aligned}
\end{equation}
recursively with initial conditions $a_{\mu,ij}=1$ and $b_{\mu,ij}=1$ for the root node. The complexity for calculating these coefficients is $\bigO(M(\bar{r}^\tau)^4)$. Let us compare this with our matrix integrator proposed in \cite{Einkemmer_2024}: The complexity of calculating the coefficients $c^0_{\mu,ij}(\vb{x}^0)$ and $d^0_{\mu,ij}(\vb{x}^0)$ was $\bigO(Mr^2n^0n^1)$. Therefore, we expect that our TTN integrator is faster for matrices than the original matrix integrator, provided that the rank is small. In Section~\ref{sec:experiments} we will see that this is indeed the case.

With the new coefficients the evolution equation~(\ref{eq:K0-TTN}) for the K step becomes
\begin{equation}\label{eq:K0-TTN-new}
    \partial_{t}K^{\tau_0}_{i_{\tau_0}}(t,\vb{x}^{\tau_0})=\sum_\mu \sum_{j_{\tau_0}}\left( \alpha^{\tau_0}_\mu(\vb{x}^{\tau_0}-\vb{\nu}^{\tau_0}_\mu) a^{\tau_0}_{\mu,i_{\tau_0} j_{\tau_0}} K_{j_{\tau_0}}(t,\vb{x}^{\tau_0}-\vb{\nu}^{\tau_0}_{\mu}) - \alpha^{\tau_0}_\mu(\vb{x}^{\tau_0}) b^{\tau_0}_{\mu,i_{\tau_0} j_{\tau_0}} K_{j_{\tau_0}}(t,\vb{x}^{\tau_0}) \right).
\end{equation}
The computational cost for the evolution equation of the K step is therefore $\bigO(M n^{\tau_0} (r^{\tau_0})^2)$.

After integrating Equation~(\ref{eq:K0-TTN-new}) and performing the QR decomposition via Equation~(\ref{eq:K0-TTN-QR}), we obtain updated low-rank factors $X^{\tau_0}_{0,i_{\tau_0}}(\vb{x}^{\tau_0})$ and we use them for updating the inner products of $A^{\tau_0}_{\mu, i_{\tau_0} j_{\tau_0}}$ and $B^{\tau_0}_{\mu, i_{\tau_0} j_{\tau_0}}$.
With the updated $A^{\tau_0}_{\mu,i_{\tau_0} j_{\tau_0}}$ and $B^{\tau_0}_{\mu,i_{\tau_0} j_{\tau_0}}$ as basic building blocks, the evolution equation~(\ref{eq:S0-TTN}) of the S step becomes
\begin{equation}\label{eq:S0-TTN-new}
    \frac{\diff{}}{\diff{t}}S^{\tau_0}_{i_{\tau_0} j_{\tau_0}}(t) = -\sum_\mu \sum_{k_{\tau_0},l_{\tau_0}} S^{\tau_0}_{k_{\tau_0} l_{\tau_0}}(t) \left( A^{\tau_0}_{\mu, i_{\tau_0} k_{\tau_0}} a^{\tau_0}_{\mu,j_{\tau_0} l_{\tau_0}} - B^{\tau_0}_{\mu, i_{\tau_0} k_{\tau_0}} b^{\tau_0}_{\mu,j_{\tau_0} l_{\tau_0}} \right).
\end{equation}
Therefore, the computational cost for the evolution equation of the S step scales with $\bigO(M(r^{\tau_0})^3)$.

For the second step of Section~\ref{sec:ttn-algorithm}, where the right child node $\tau_1$ is updated, the equations are modified in a similar fashion.

For the third step, where $C^\tau_{i_\tau i_{\tau_0} i_{\tau_1} j_\tau j_{\tau_0} j_{\tau_1}}$ is updated, the evolution equation (\ref{eq:C-TTN}) can be written as
\begin{equation}\label{eq:C-TTN-new}
    \frac{\diff{}}{\diff{t}}C^{\tau}_{i_\tau i_{\tau_0} i_{\tau_1}}(t) = \sum_\mu \sum_{j_\tau} \sum_{j_{\tau_0}} \sum_{j_{\tau_1}} C^{\tau}_{j_\tau j_{\tau_0} j_{\tau_1}}(t) \left( A^{\tau_0}_{\mu,i_{\tau_0} j_{\tau_0}} A^{\tau_1}_{\mu,i_{\tau_1} j_{\tau_1}} a^\tau_{\mu,i_\tau j_\tau} - B^{\tau_0}_{\mu,i_{\tau_0} j_{\tau_0}} B^{\tau_1}_{\mu,i_{\tau_1} j_{\tau_1}} b^\tau_{\mu,i_\tau j_\tau} \right).
\end{equation}
The complexity for updating the connection tensor is therefore $\bigO(M(\bar{r}^\tau)^4)$.


Besides the factorization property we also exploit in our implementation the fact that most reactions only depend on a small subset of all species. Due to this behaviour the propensity functions only depend on the population numbers of the species which are participating in the corresponding reaction. We can thus directly store the relatively few values of the propensities at the leaves of the tree and thereby avoid expensive function evaluations for every time step.

We summarize the resulting first-order integration scheme in Algorithms~\ref{alg:TTN-integrator}--\ref{alg:Subflow-Psi}: The subflows $\Phi^0(\tau)$ and $\Phi^1(\tau)$ (Algorithms~\ref{alg:Subflow-Phi-0} and \ref{alg:Subflow-Phi-1}) correspond to the first and second step described in Section~\ref{sec:ttn-algorithm} and update the low-rank factors of the child nodes $\tau_0$ and $\tau_1$. The subflow $\Psi(\tau)$ (Algorithm~\ref{alg:Subflow-Psi}) updates the tensor $C^{\tau}_{i_\tau i_{\tau_0} i_{\tau_1}}$. Algorithm~\ref{alg:TTN-integrator} computes the composition $\Psi(\tau) \circ \Phi^1(\tau) \circ \Phi^0(\tau)$ of the three subflows for node $\tau$, and we call it at the root to obtain $P(\Delta t, \vb{x}) \approx \sum_i X_{1,i}(\vb{x})$ at time $\Delta t$.

\begin{algorithm}[H]
    \caption{First-order projector-splitting tree tensor network (PS-TTN) integrator for node $\tau$\label{alg:TTN-integrator}}
    \begin{tabular}{ll}
        \textbf{Input:} & Precomputed coefficients~(\ref{eq:A_coeff}),\\
        & $X^{\tau}_{0,i_\tau}(\vb{x}^\tau)=\sum_{i_{\tau_0}} \sum_{i_{\tau_1}}  Q^\tau_{0,i_\tau i_{\tau_0} i_{\tau_1}} X^{\tau_0}_{0,i_{\tau_0}}(\vb{x}^{\tau_0}) X^{\tau_1}_{0,i_{\tau_1}}(\vb{x}^{\tau_1})$ and,\\
        & for the root node, $C_{0,i i_0 i_1} = Q_{0,i i_0 i_1}$\\

        \textbf{Output:} & $X^{\tau}_{1,i_\tau}(\vb{x}^\tau)=\sum_{i_{\tau_0}} \sum_{i_{\tau_1}}  Q^\tau_{1,i_\tau i_{\tau_0} i_{\tau_1}} X^{\tau_0}_{1,i_{\tau_0}}(\vb{x}^{\tau_0}) X^{\tau_1}_{1,i_{\tau_1}}(\vb{x}^{\tau_1})$\\
    \end{tabular}

    \begin{algorithmic}[1]
        \State Compute \Call{$\Phi^0$}{$\tau$}
        \State Compute \Call{$\Phi^1$}{$\tau$}
        \State Compute \Call{$\Psi$}{$\tau$}
    \end{algorithmic}
\end{algorithm}

\begin{algorithm}[H]
    \caption{Subflow $\Phi^{0}(\tau)$\label{alg:Subflow-Phi-0}}
    \begin{tabular}{ll}
        \textbf{Input:} & $C^\tau_{0,i_\tau i_{\tau_0} i_{\tau_1}}$, $X^{\tau_0}_{0,i_{\tau_0}}(\vb{x}^{\tau_0})$ and $X^{\tau_1}_{0,i_{\tau_1}}(\vb{x}^{\tau_1})$\\

        \textbf{Output:} & $C^{\tau_0}_{0,i_{\tau_0} i_{\tau_{00}} i_{\tau_{01}}}$, $C^\tau_{1,i_\tau i_{\tau_0} i_{\tau_1}}$ and $X^{\tau_0}_{1,i_{\tau_0}}(\vb{x}^{\tau_0})$\\
    \end{tabular}

    \begin{algorithmic}[1]
        \State Decompose $C^\tau_{0,i_\tau i_{\tau_0} i_{\tau_1}} = \sum_{j_{\tau_0}} G^{\tau_0}_{i_\tau j_{\tau_0} i_{\tau_1}} S^{\tau_0}_{0,i_{\tau_0} j_{\tau_0}}$ via a QR factorization
        \If{$\tau_0$ is an internal node}
        \State Calculate $C^{\tau_0}_{0,i_{\tau_0} i_{\tau_{00}} i_{\tau_{01}}}=\sum_{j_{\tau_0}} Q^{\tau_0}_{0,j_{\tau_0} i_{\tau_{00}} i_{\tau_{01}}} S^{\tau_0}_{0,j_{\tau_0} i_{\tau_0}}$
        \State Call \Call{PS-TTN integrator}{$\tau_0$} using Algorithm~\ref{alg:TTN-integrator}
        \State Decompose $C^{\tau_0}_{3,i_{\tau_0} i_{\tau_{00}} i_{\tau_{01}}} = \sum_{j_{\tau_0}} Q^{\tau_0}_{1,j_{\tau_0} i_{\tau_{00}} i_{\tau_{01}}} \tilde{S}^{\tau_0}_{j_{\tau_0} i_{\tau_0}}$ via a QR factorization
        \State Update $A^{\tau_0}_{\mu,i_{\tau_0} j_{\tau_0}}$ and $B^{\tau_0}_{\mu,i_{\tau_0} j_{\tau_0}}$ with $X^{\tau_0}_{1,i_{\tau_0}}(\vb{x}^{\tau_0})$ using Equation~(\ref{eq:A_coeff-recursion})
        \Else
        \State Calculate $a^{\tau_0}_{\mu,i_{\tau_0} j_{\tau_0}}$ and $b^{\tau_0}_{\mu,i_{\tau_0} j_{\tau_0}}$ using Equation~(\ref{eq:K0-TTN_coeff-new})
        \State Integrate $K^{\tau_0}_{i_{\tau_0}}$ from $0$ to $\Delta t$ with initial value $K^{\tau_0}_{i_{\tau_0}}(0,\vb{x}^{\tau_0}) = \sum_{j_{\tau_0}} X^{\tau_0}_{0,j_{\tau_0}}(\vb{x}^{\tau_0}) S^{\tau_0}_{0,j_{\tau_0}i_{\tau_0}}$ using Equation~(\ref{eq:K0-TTN-new})
        \State Decompose $K^{\tau_0}_{i_{\tau_0}}(\Delta t,\vb{x}^{\tau_0}) = \sum_{j_{\tau_0}} X^{\tau_0}_{1,j_{\tau_0}}(\vb{x}^{\tau_0}) \tilde{S}^{\tau_0}_{j_{\tau_0}i_{\tau_0}}$ via a QR factorization
        \State Update $A^{\tau_0}_{\mu,i_{\tau_0} j_{\tau_0}}$ and $B^{\tau_0}_{\mu,i_{\tau_0} j_{\tau_0}}$ with $X^{\tau_0}_{1,i_{\tau_0}}(\vb{x}^{\tau_0})$ using Equation~(\ref{eq:A_coeff})
        \EndIf
        \State Integrate $S^{\tau_0}_{i_{\tau_0} j_{\tau_0}}$ from $0$ to $\Delta t$ with initial value $S_{i_{\tau_0} j_{\tau_0}}(0) = \tilde{S}^{\tau_0}_{i_{\tau_0} j_{\tau_0}}$ using Equation~(\ref{eq:S0-TTN-new}) and set $S^{\tau_0}_{1,i_{\tau_0} j_{\tau_0}} = S^{\tau_0}_{i_{\tau_0} j_{\tau_0}}(\Delta t)$
        \State Calculate $C^\tau_{1,i_{\tau} i_{\tau_0} i_{\tau_1}} = \sum_{j_{\tau_0}} G^{\tau_0}_{i_\tau j_{\tau_0} i_{\tau_1}} S^{\tau_0}_{1,i_{\tau_0} j_{\tau_0}}$
    \end{algorithmic}
\end{algorithm}

\begin{algorithm}[H]
    \caption{Subflow $\Phi^{1}(\tau)$\label{alg:Subflow-Phi-1}}
    \begin{tabular}{ll}
        \textbf{Input:} & $C^\tau_{0,i_\tau i_{\tau_0} i_{\tau_1}}$, $X^{\tau_0}_{1,i_{\tau_0}}(\vb{x}^{\tau_0})$ and $X^{\tau_0}_{0,i_{\tau_0}}(\vb{x}^{\tau_0})$\\

        \textbf{Output:} & $C^{\tau_1}_{0,i_{\tau_1} i_{\tau_{10}} i_{\tau_{11}}}$, $C^\tau_{2,i_\tau i_{\tau_0} i_{\tau_1}}$ and $X^{\tau_1}_{1,i_{\tau_1}}(\vb{x}^{\tau_1})$\\
    \end{tabular}

    \begin{algorithmic}[1]
        \State Decompose $C^\tau_{1,i_\tau i_{\tau_0} i_{\tau_1}} = \sum_{j_{\tau_1}} G^{\tau_1}_{i_\tau i_{\tau_0} j_{\tau_1}} S^{\tau_1}_{0,i_{\tau_1} j_{\tau_1}}$ via a QR factorization
        \If{$\tau_1$ is an internal node}
        \State Calculate $C^{\tau_1}_{0,i_{\tau_1} i_{\tau_{10}} i_{\tau_{11}}}=\sum_{j_{\tau_1}} Q^{\tau_1}_{0,j_{\tau_1} i_{\tau_{10}} i_{\tau_{11}}} S^{\tau_1}_{0,j_{\tau_1} i_{\tau_1}}$
        \State Call \Call{PS-TTN integrator}{$\tau_1$} 
        \State Decompose $C^{\tau_1}_{3,i_{\tau_1} i_{\tau_{10}} i_{\tau_{11}}} = \sum_{j_{\tau_1}} Q^{\tau_1}_{1,j_{\tau_1} i_{\tau_{10}} i_{\tau_{11}}} \tilde{S}^{\tau_1}_{j_{\tau_1} i_{\tau_1}}$ via a QR factorization
        \State Update $A^{\tau_1}_{\mu,i_{\tau_1} j_{\tau_1}}$ and $B^{\tau_1}_{\mu,i_{\tau_1} j_{\tau_1}}$ with $X^{\tau_1}_{1,i_{\tau_1}}(\vb{x}^{\tau_1})$ 
        \Else
        \State Calculate $a^{\tau_1}_{\mu,i_{\tau_1} j_{\tau_1}}$ and $b^{\tau_1}_{\mu,i_{\tau_1} j_{\tau_1}}$ 
        \State Integrate $K^{\tau_1}_{i_{\tau_1}}$ from $0$ to $\Delta t$ with initial value $K^{\tau_1}_{i_{\tau_1}}(0,\vb{x}^{\tau_1}) = \sum_{j_{\tau_1}} X^{\tau_1}_{0,j_{\tau_1}}(\vb{x}^{\tau_1}) S^{\tau_1}_{0,j_{\tau_1}i_{\tau_1}}$ 
        \State Decompose $K^{\tau_1}_{i_{\tau_1}}(\Delta t,\vb{x}^{\tau_1}) = \sum_{j_{\tau_1}} X^{\tau_1}_{1,j_{\tau_1}}(\vb{x}^{\tau_1}) \tilde{S}^{\tau_1}_{j_{\tau_1}i_{\tau_1}}$ via a QR factorization
        \State Update $A^{\tau_1}_{\mu,i_{\tau_1} j_{\tau_1}}$ and $B^{\tau_1}_{\mu,i_{\tau_1} j_{\tau_1}}$ with $X^{\tau_1}_{1,i_{\tau_1}}(\vb{x}^{\tau_1})$ 
        \EndIf
        \State Integrate $S^{\tau_1}_{i_{\tau_1} j_{\tau_1}}$ from $0$ to $\Delta t$ with initial value $S_{i_{\tau_1}j_{\tau_1}}(0) = \tilde{S}^{\tau_1}_{i_{\tau_1} j_{\tau_1}}$ and set $S^{\tau_1}_{1,i_{\tau_1} j_{\tau_1}} = S^{\tau_1}_{i_{\tau_1} j_{\tau_1}}(\Delta t)$ 
        \State Calculate $C^\tau_{2,i_{\tau} i_{\tau_1} i_{\tau_1}} = \sum_{j_{\tau_1}} G^{\tau_1}_{i_\tau i_{\tau_0} j_{\tau_1}} S^{\tau_1}_{1,i_{\tau_1} j_{\tau_1}}$
    \end{algorithmic}
\end{algorithm}

\begin{algorithm}[H]
    \caption{Subflow $\Psi(\tau)$\label{alg:Subflow-Psi}}
    \begin{tabular}{ll}
        \textbf{Input:} & $C^\tau_{2,i_\tau i_{\tau_0} i_{\tau_1}}$, $X^{\tau_0}_{1,i_{\tau_0}}(\vb{x}^{\tau_0})$ and $X^{\tau_0}_{1,i_{\tau_0}}(\vb{x}^{\tau_0})$\\

        \textbf{Output:} & $C^\tau_{3,i_\tau i_{\tau_0} i_{\tau_1}}$\\
    \end{tabular}

    \begin{algorithmic}[1]
    \State Integrate $C^\tau_{i_\tau i_{\tau_0} i_{\tau_1}}$ from $0$ to $\Delta t$ with initial value $C^\tau_{i_\tau i_{\tau_0} i_{\tau_1}}(0) = C^\tau_{2,i_\tau i_{\tau_0} i_{\tau_1}}$ using Equation~(\ref{eq:C-TTN-new}) \State Set $C^\tau_{3,i_\tau i_{\tau_0} i_{\tau_1}} = C^\tau_{i_\tau i_{\tau_0} i_{\tau_1}}(\Delta t)$ (or $Q_{i i_0 i_1} = C_{i i_0 i_1}(\Delta t)$, if $\tau$ is the root)
    \end{algorithmic}
\end{algorithm}

Since reaction networks often include reactions with different time scales, the resulting equations can become stiff. Our \CC{} implementation gives the user therefore the flexibility to choose between explicit and implicit Euler schemes for the time integration. The implicit scheme uses the matrix-free GMRES implementation of \texttt{Eigen} \cite{Saad_1986,Guennebaud_2010}. We also use \texttt{Ensign}~\cite{Cassini_2021} for the low-rank data structures and other routines related to the dynamical low-rank approximation.


\section{Numerical experiments}\label{sec:experiments}
We test our implementation with two models from the field of biochemistry. The smaller model, the bacteriophage-$\lambda$ (“lambda phage”), was primarily chosen to validate the implementation and to investigate the approximation accuracy. The relatively small system size of this model allows for computing a reference solution of the CME without a low-rank approximation, which we compare with our TTN integrator. The second example, the BAX pore assembly model, is a large reaction network and computing a reference solution of the CME without any approximation is no longer feasible. We therefore compare this model with a state-of-the-art implementation of SSA called \texttt{GillesPy2} \cite{Matthew_2023}.



\subsection{Lambda phage}
In this first example, we apply our TTN integrator to the model for the life cycle of the lambda phage as described in \cite{Hegland_2007}. This model consists of five different species, which can interact by ten different reactions (Table~\ref{tab:lp}). Recall that species which are of no further interest are denoted by $\star$.

\begin{table}[H]
    \centering
    \begin{tabular}{c|c|c}
        \hline 
        $\mu$ & Reaction $R_\mu$ & Propensity function $\alpha_\mu(\vb{x})$\tabularnewline
        \hline 
        0 & \ce{$\star$ -> $S_{0}$} & $a_{0}b_{0}/(b_{0}+x_{1})$\tabularnewline
        1 & \ce{$\star$ -> $S_{1}$} & $(a_{1}+x_{4})b_{1}/(b_{1}+x_{0})$\tabularnewline
        2 & \ce{$\star$ -> $S_{2}$} & $a_{2}b_{2}x_{1}/(b_{2}x_{1}+1)$\tabularnewline
        3 & \ce{$\star$ -> $S_{3}$} & $a_{3}b_{3}x_{2}/(b_{3}x_{2}+1)$\tabularnewline
        4 & \ce{$\star$ -> $S_{4}$} & $a_{4}b_{4}x_{2}/(b_{4}x_{2}+1)$\tabularnewline
        5 & \ce{$S_{0}$ -> $\star$} & $c_{0}\cdot x_{0}$\tabularnewline
        6 & \ce{$S_{1}$ -> $\star$} & $c_{1}\cdot x_{1}$\tabularnewline
        7 & \ce{$S_{2}$ -> $\star$} & $c_{2}\cdot x_{2}$\tabularnewline
        8 & \ce{$S_{3}$ -> $\star$} & $c_{3}\cdot x_{3}$\tabularnewline
        9 & \ce{$S_{4}$ -> $\star$} & $c_{4}\cdot x_{4}$\tabularnewline
        \hline 
        \end{tabular}\qquad{}%
        \begin{tabular}{c|c|c|c|c|c}
        \hline 
         & $i=0$ & $i=1$ & $i=2$ & $i=3$ & $i=4$\tabularnewline
        \hline 
        $a_i$ & $0.5$ & $1$ & $0.15$ & $0.3$ & $0.3$\tabularnewline
        $b_i$ & $0.12$ & $0.6$ & $1$ & $1$ & $1$\tabularnewline
        $c_i$ & $0.0025$ & $0.0007$ & $0.0231$ & $0.01$ & $0.01$\tabularnewline
        \hline 
    \end{tabular}
    \caption{Reactions, propensity functions and parameters of the lambda phage
    model.\label{tab:lp}}
\end{table}

The life cycle of the lambda phage represents a naturally occurring toggle switch, which is the analogue of the flip-flop circuit in electronics (for more details on the toggle switch, we refer to \cite{Gardner_2000}). The lambda phage infects \emph{E. coli}, and depending on the environment, either stays dormant in the bacterial host (\emph{lysogenic phase}) or multiplies, reassembles itself and breaks out of the host (\emph{lytic phase}). If enough $S_{4}$ is present in the environment, $S_{1}$ is produced and the system is in the lysogenic phase. Abundance of $S_{1}$ in turn inhibits the formation of $S_{0}$ via reaction~$R_0$. If the amount of $S_{4}$ in the environment is scarce, the production of $S_{0}$ causes the system to enter the lytic phase and the transcription of new copies of $S_{1}$ via reaction~$R_1$ is inhibited.

As an initial value the multinomial distribution with parameters $n=3$ and $p=(0.05,\dots,0.05)$ has been chosen:
\begin{equation*}
    P(0,\vb{x})=
    \begin{cases}
        \frac{3!}{x_{0}!\cdots x_{4}!(3-|\vb{x}|)!}0.05^{|\vb{x}|}(1-5\cdot0.05)^{3-|\vb{x}|} & \text{if}\quad|\vb{x}|\le3,\\
        0 & \text{else,}
    \end{cases}
\end{equation*}
where $|\vb{x}|=x_0+\dots+x_4$. We solved the CME on the time interval $[0,10]$ with truncation indices $\vb{\eta}=(0,0,0,0,0)$ and $\vb{\zeta}=(15,40,10,10,10)$.

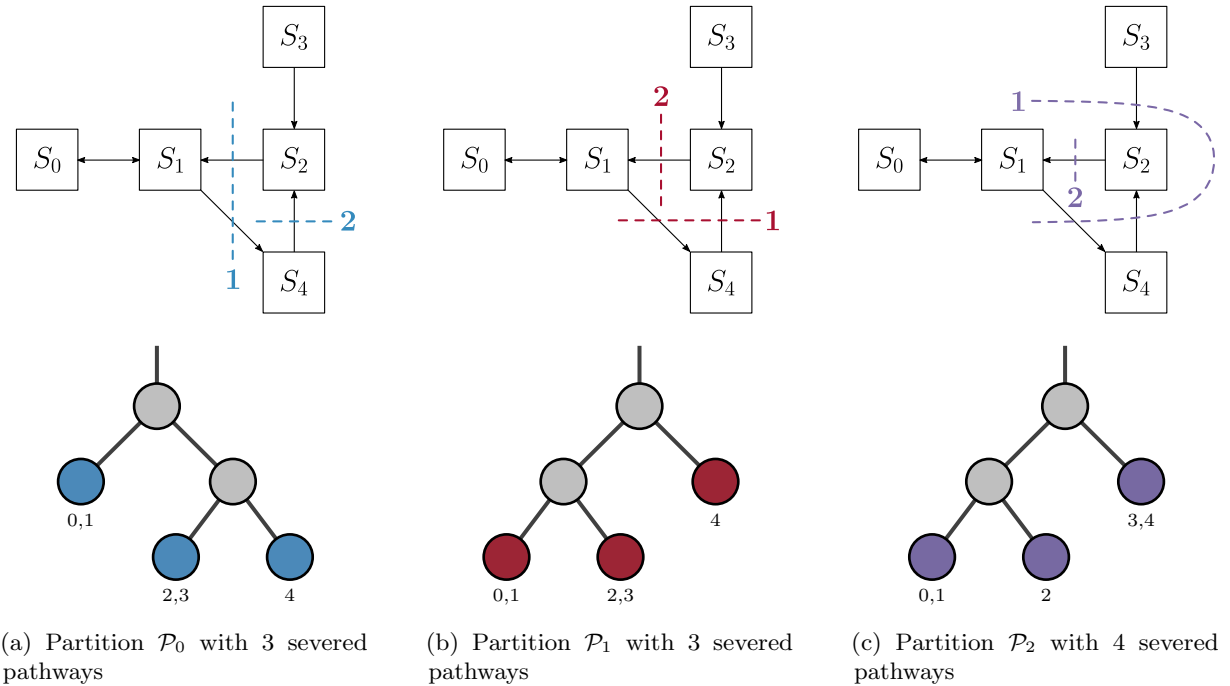
\begin{figure}[H]
    \centering
    \begin{subfigure}[b]{0.3\textwidth} 
        \centering
        \includegraphics[height=0.85\textwidth]{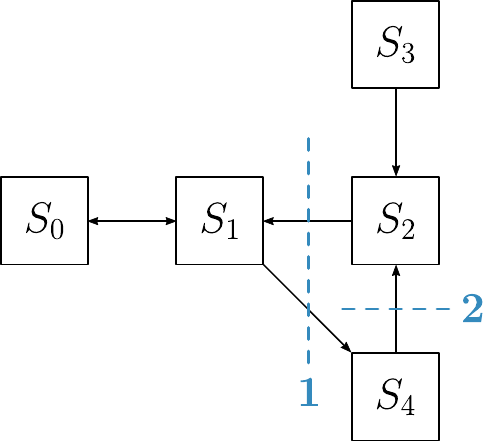}
        \\[1ex]
        \begin{tikzpicture}       
    \InternalNode[]{root}
    \node[above=0.5cm of root](D){};
    \ExternalNode[x=-1,y=-1,label={0,1},position=below]{0}
    \InternalNode[x=1,y=-1]{1}
    \ExternalNode[x=0.25,y=-2,label={2,3},position=below]{10}
    \ExternalNode[x=1.75,y=-2,label={4},position=below]{11}
    
    \Edge(D)(root)
    \Edge(root)(0)
    \Edge(root)(1)
    \Edge(1)(10)
    \Edge(1)(11)
\end{tikzpicture}
        \caption{Partition $\mathcal{P}_0$ with 3 severed pathways\label{fig:lp-partition-0}}
    \end{subfigure}
    \hfill
    \begin{subfigure}[b]{0.3\textwidth}
        \centering
        \includegraphics[height=0.85\textwidth]{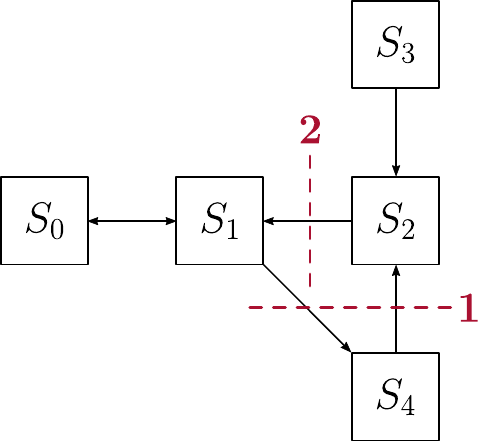}
        \\[1ex]
        \begin{tikzpicture}       
    \InternalNode[]{root}
    \node[above=0.5cm of root](D){};
    \InternalNode[x=-1,y=-1,]{0}
    \Vertex[x=1,y=-1,label={4},position=below,color=custom_red]{1}
    \Vertex[x=-1.75,y=-2,label={0,1},position=below,color=custom_red]{00}
    \Vertex[x=-0.25,y=-2,label={2,3},position=below,color=custom_red]{01}
    
    \Edge(D)(root)
    \Edge(root)(0)
    \Edge(root)(1)
    \Edge(0)(00)
    \Edge(0)(01)
\end{tikzpicture}
        \caption{Partition $\mathcal{P}_1$ with 3 severed pathways}
    \end{subfigure}
    \hfill
    \begin{subfigure}[b]{0.3\textwidth}
        \centering
        \includegraphics[height=0.85\textwidth]{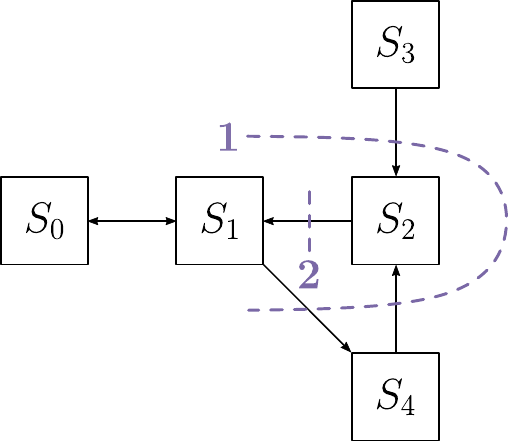}
        \\[1ex]
        \begin{tikzpicture}       
    \InternalNode[]{root}
    \node[above=0.5cm of root](D){};
    \InternalNode[x=-1,y=-1,]{0}
    \Vertex[x=1,y=-1,label={3,4},position=below,color=custom_violet]{1}
    \Vertex[x=-1.75,y=-2,label={0,1},position=below,color=custom_violet]{00}
    \Vertex[x=-0.25,y=-2,label={2},position=below,color=custom_violet]{01}
    
    \Edge(D)(root)
    \Edge(root)(0)
    \Edge(root)(1)
    \Edge(0)(00)
    \Edge(0)(01)
\end{tikzpicture}
        \caption{Partition $\mathcal{P}_2$ with 4 severed pathways}
    \end{subfigure}
    \caption{Dependency graphs (top) and tree structures (bottom) for the three partitions of the lambda phage example. The notation $A\longrightarrow B$ in the dependency graphs indicates that species $A$ depends on species $B$ due to a specific reaction propensity. Small numbers at the leaves of the tree structure denote on which species the corresponding low-rank factors depend.\label{fig:lp-partitions}}
\end{figure}

\subsubsection{Benchmark results}
For our benchmark we compute solutions with partition $\mathcal{P}_0$ shown in Figure~\ref{fig:lp-partition-0}. For the time integration we use, unless otherwise mentioned, an explicit Euler scheme with time step size $\Delta t=10^{-3}$. All computations were performed on a MacBook Pro with a $2\,\mathrm{GHz}$ Intel Core i5 Skylake (6360U) processor.
We compare our solution with an ``exact'' reference solution, which was obtained by solving the full CME on the truncated state space with \texttt{scipy.solve\_ivp}. Due to the relatively large system size the computation of the full solution is very costly and therefore substantially slower than the DLR approximation. We also compare our method with SSA for $10^4$, $10^5$ and $10^6$ generated trajectories or runs.

Note that for this example we require $r_{\tau_0}=r_{\tau_1}$ for both subtrees $\tau_0$ and $\tau_1$ of a given tree $\tau$ (this is a special case of the rank condition (\ref{eq:Tucker-condition})). This means that we assign a single rank for every cut of the reaction system and the two resulting partitions are approximated to the same degree. With this convention we can concisely write the ranks as an ordered tuple, which has the same ordering as the cuts. Let us consider partition $\mathcal{P}_0$ (shown in Figure~\ref{fig:lp-partition-0}) as an example: Here we write the rank as $r=(r_\text{cut 1},r_\text{cut 2})$, with $r_\text{cut 1}=r^0=r^1$ and $r_\text{cut 2}=r^{10}=r^{11}$.

\begin{figure}[H]
    \centering
    \includegraphics[scale=0.8]{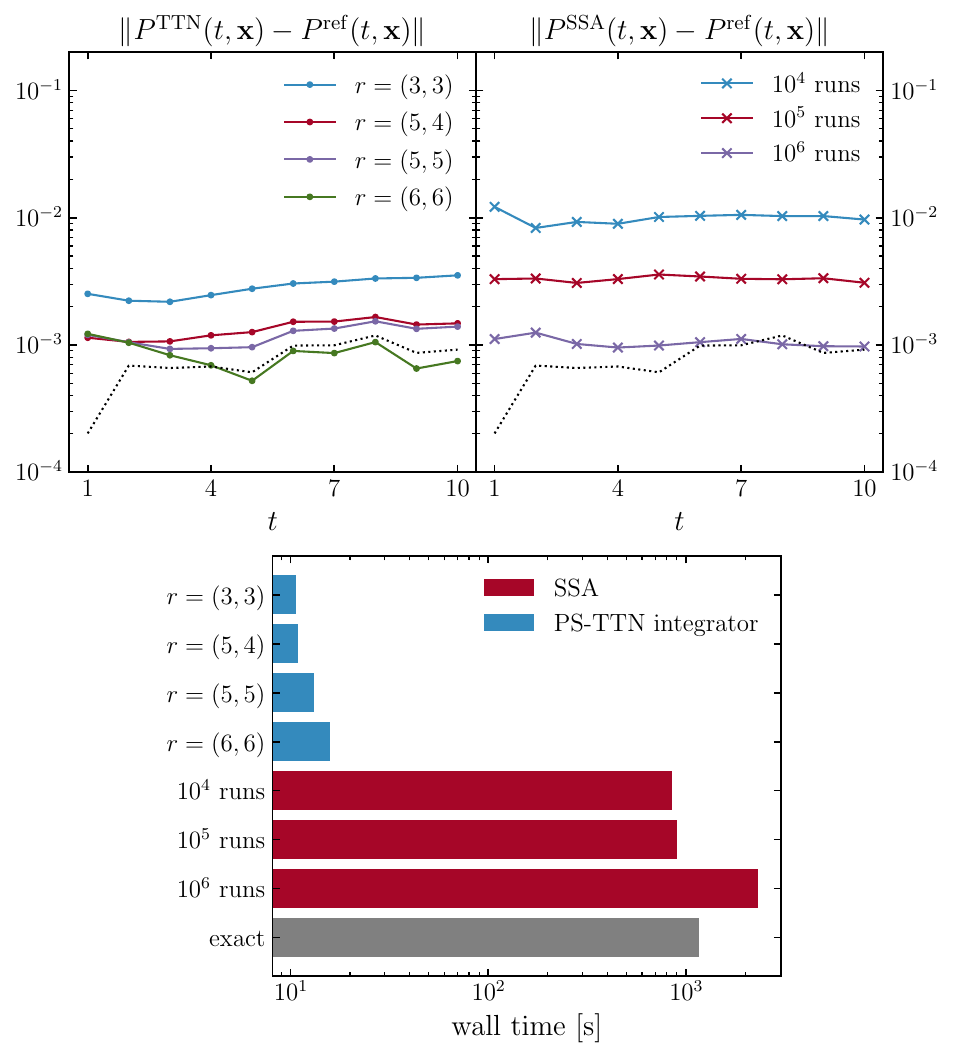}
    \caption{Top: Time-dependent 2-norm error for the PS-TTN integrator (left panel) and SSA (right panel) for the lambda phage model. The reference solution was obtained by solving the full CME on the truncated state space with \texttt{scipy.solve\_ivp}. For the PS-TTN we used explicit Euler with time step size $\Delta t=10^{-3}$. For convenience we show the best-approximation in the matrix case for a truncation at $r=5$, which was obtained by a truncated SVD of the reference solution (dotted line). Bottom: Wall time comparison (in seconds) for the PS-TTN integrator, SSA and the reference solution (``exact'').\label{fig:lp-ttn-ssa}}
\end{figure}

Since a reference solution is at our disposal, we can compare our PS-TTN integrator directly with SSA in terms of accuracy and performance. In the top panels of Figure~\ref{fig:lp-ttn-ssa}, we plot the time-dependent absolute error (with respect to the 2-norm) of our PS-TTN integrator for different ranks and for SSA for different numbers of trajectories. The theoretically best result of the DLR approximation is the truncated SVD of the reference solution. For convenience we show this best-approximation in the matrix case for a truncation at $r=5$.
Already for relatively small ranks $r=(5,4)$ and $r=(5,5)$ our solutions are as accurate as SSA with $10^6$ runs, and the solution with $r=(6,6)$ comfortably outperforms the Monte Carlo method. The bottom panel of Figure~\ref{fig:lp-ttn-ssa} sheds a light on the wall time (on a logarithmic scale) of the different methods. The PS-TTN integrator is for all settings almost two magnitudes faster than the integrator for the full CME and the SSA. For $r=(5,5)$, the wall time of approximately $10\,\mathrm{s}$ is also by a factor of $5$ faster than our matrix integrator as proposed in \cite{Einkemmer_2024}, with a wall time of $53\,\mathrm{s}$ for $r=5$. This comes from the fact that the matrix integrator does not exploit the factorization property.

\begin{figure}[H]
    \centering
    \includegraphics[scale=0.8]{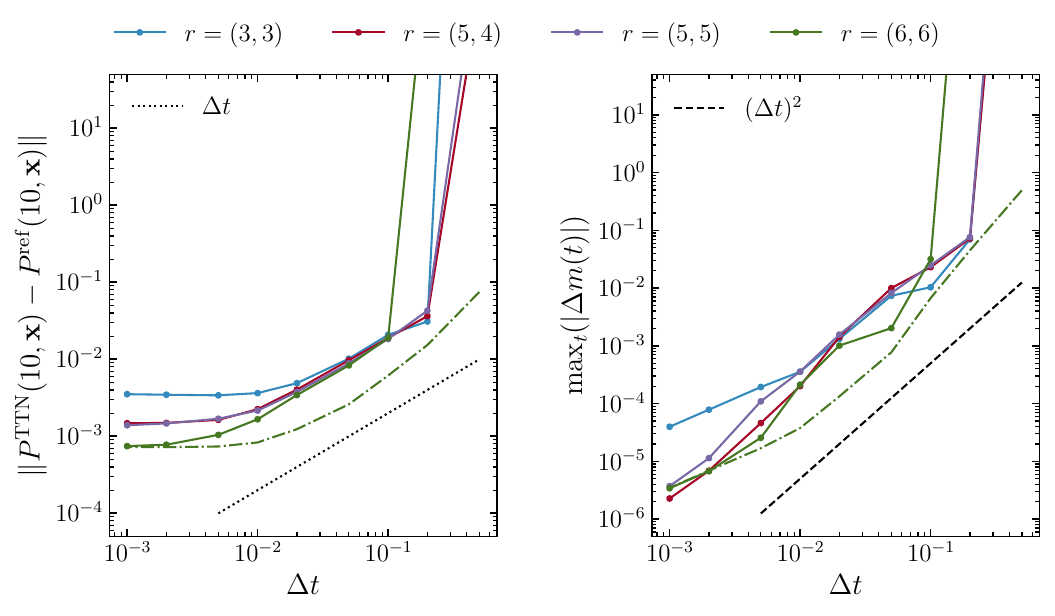}
    \caption{Left: 2-norm error for the PS-TTN solution at time $t=10$ depending on the time step size $\Delta t$ for the lambda phage model. For the solid lines, an explicit Euler scheme was used for the time integration. Right: Maximum mass error over the time integration interval $[0, 10]$ depending on the time step size $\Delta t$. For sake of comparison we show also a solution for the PS-TTN integrator with an implicit Euler scheme and rank $r=(6,6)$ (green, dash-dotted line).\label{fig:lp-time-err}}
\end{figure}

In Section~\ref{sec:ttn}, we mentioned that our PS-TTN integrator is first-order accurate in time. This is indeed the case for a certain time step size range, as can be seen in the left panel of Figure~\ref{fig:lp-time-err}, where we show the error between our PS-TTN solution and the reference solution depending on the time step size $\Delta t$ at time $t=10$. If the time step size is small, then then the error of the low-rank approximation becomes dominant, and for higher ranks (e.g., $r=(6,6)$) this effect can be seen only for smaller time step sizes ($\Delta t=10^{-2}$). Recall that the PS-TTN solutions shown in Figure~\ref{fig:lp-ttn-ssa} were computed with $\Delta t=10^{-3}$, thus we directly compare the error of the low-rank approximation with the $1/\sqrt{N}$-behaviour of the Monte Carlo method.

For rank $r=(6,6)$, we observe instabilities for smaller time step sizes. The origin of this behaviour lies in our explicit Euler scheme. For a time integration with implicit Euler, this effect is indeed mitigated (green, dash-dotted line).

Since we do not explicitly conserve mass in the sense of Equation~(\ref{eq:mass}) with our numerical scheme, it is also important to study the mass error of our method. In the right panel of Figure~\ref{fig:lp-time-err} we plot the maximum mass error over the time integration interval $[0,10]$ depending on the time step size $\Delta t$. We can see that in the stable regime the mass error is more or less rank-independent and it seems to scale with $\bigO(\Delta t^{2})$. For time step size $\Delta t=10^{-3}$ the maximum mass error is smaller than $10^{-5}$ and negligible compared to the low-rank approximation error.

Finally, we report the memory requirements for storing the TTN representation of the probability distribution at a single point in time. For the full probability distribution $6.99\,\mathrm{MB}$ of storage is needed, which is reduced to $32.7\,\mathrm{kB}$ ($\approx 0.47\%$ of the original requirements) for $r=(5,5)$ and $39.8\,\mathrm{kB}$ ($\approx 0.57\%$ of the original requirements) for $r=(6,6)$.

\subsubsection{Partition study}
So far we only reported benchmarks for the lambda phage model with partition $\mathcal{P}_0$. However, it is expected that the accuracy (and for larger problems, also the wall time and the memory requirements) are related to the cuts of the reaction pathways. Therefore, we will compare in this section the partition $\mathcal{P}_0$ with the two other partitions shown in Figure~\ref{fig:lp-partitions}. The three different partitions, which result from cutting the reaction network twice, are depicted schematically in Figure~\ref{fig:lp-partitions}.

\begin{figure}[H]
    \centering
    \includegraphics[scale=0.8]{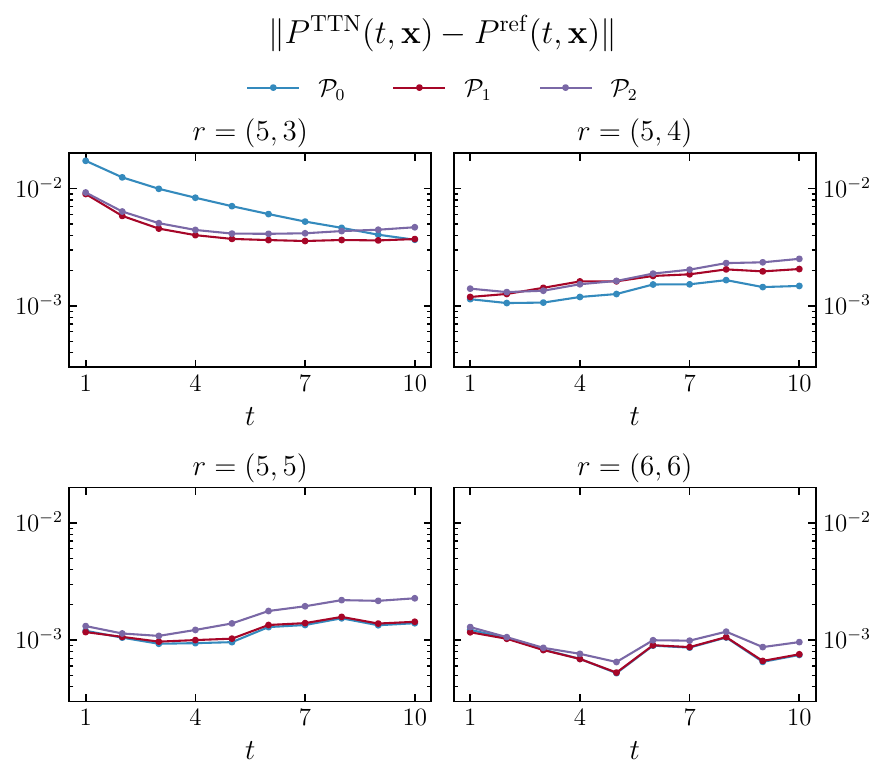}
    \caption{Time-dependent 2-norm error for the PS-TTN integrator using explicit Euler with time step size $\Delta t=10^{-3}$ for the three different partitions of the lambda phage model depicted in Figure~\ref{fig:lp-partitions}.\label{fig:lp-partition-comparison}}
\end{figure}

First we will address the influence of the chosen partition on the accuracy of our numerical scheme. To this end, we show the time-dependent absolute error between the PS-TTN integrator using explicit Euler with time step size $\Delta t=10^{-3}$ for the three different partitions with different ranks in Figure~\ref{fig:lp-partition-comparison}. Looking at the number of cuts shown in Figure~\ref{fig:lp-partitions}, we expect that partition $\mathcal{P}_2$ needs a higher rank to achieve a similar accuracy as the other two partitions. This is indeed the case for $r=(5,4)$ and $r=(5,5)$. Surprisingly, for $r=(5,3)$ partition $\mathcal{P}_0$ shows (up to $t=8$) the largest error. If we take a look at Figure~\ref{fig:lp-partitions}, we see that the cut reaction pathways and the number of approximated pathways are the same for $\mathcal{P}_0$ and $\mathcal{P}_1$ (this is also the reason why both show the same accuracy for $r=(5,5)$). In the case of $r=(5,3)$ however, the pathway from $S_4$ to $S_2$ is approximated with only three low-rank factors for $\mathcal{P}_0$, whereas for $\mathcal{P}_1$ three low-rank factors are used to approximate the pathway from $S_4$ to $S_2$. By comparing the related propensity functions ($\alpha_4(\vb{x})$ for $\mathcal{P}_0$ and $\alpha_2(\vb{x})$ for $\mathcal{P}_1$), we see that for $x_2=x_1$, the propensity function $\alpha_4(\vb{x})$ is twice as large as $\alpha_2(\vb{x})$ and thus the probability for reaction $R_4$ is twice the probability for $R_2$. Apparently three low-rank factors are not enough for approximating reaction~$R_4$, whereas it is sufficient for reaction~$R_2$, which is less likely to occur. Thus the difference in the importance of the reactions is the reason why we observe an increased error for $\mathcal{P}_0$.

\subsection{Reaction cascade}\label{sec:cascade}
The second and more challenging example is a reaction system with 40 reactions and 20 species which was taken from \cite{Dolgov_2015}. The reactions are shown in Table~\ref{tab:cascade}.

\begin{table}[H]
    \centering
    \begin{tabular}{c|c|cc}
        \hline 
        $\mu$ & Reaction $R_\mu$ & Propensity function $\alpha_\mu(\vb{x})$ & 
        \tabularnewline
        \hline 
        0 & \ce{$\star$ -> $S_{0}$} & $a$ &
        \tabularnewline
        1--19 & \ce{$\star$ -> $S_{i}$} & $x_{i-1} / (b + x_{i-1})$ & $(i=1,\dots,19)$
        \tabularnewline
        20--39 & \ce{$S_{j}$ -> $\star$} & $c \cdot x_{j}$ & $(j=0,\dots,19)$
        \tabularnewline
        \hline 
    \end{tabular}
    \caption{Reactions and propensity functions of the reaction cascade system with parameters $a=0.7$, $b=5$, $c=0.07$.\label{tab:cascade}}
\end{table}

Species $S_0$ is constantly produced via reaction $R_0$ and the population number of $S_0$ directly influences the creation of $S_1$ from reaction $R_1$. More generally, the creation of species $S_i$ via $R_i$ depends explicitly on the population number of the preceding species $S_{i-1}$ for $i=1,\dots,19$, and due to this cascade-like structure this is an example of a so-called \emph{reaction cascade}. Moreover, the species in our reaction system can also decay via reactions $R_{20}$--$R_{40}$.

In our numerical experiments we used a deterministic initial condition where all population numbers are zero,
\begin{equation*}
    P(0,\vb{x})=\prod_{i=0}^{19} \delta_{x_i 0},
\end{equation*}
where $\delta_{x_i 0}$ is the Kronecker delta. The CME was solved on the time interval $[0,350]$ on the truncated state space with truncation indices $\eta_i=0$ and $\zeta_i=63$ for $i=0,\dots,19$.

\begin{figure}[H]
    \begin{subfigure}[b]{\textwidth}
        \centering
        \begin{tikzpicture}
    \InternalNode{root}
    \node[left=0.5cm of root](D){};
    \InternalNode[x=1]{1}
    \InternalNode[x=2]{11}
    \InternalNode[x=3]{111}
    \InternalNode[x=4]{1111}
    \InternalNode[x=5]{11111}
    \InternalNode[x=6]{111111}
    \InternalNode[x=7]{1111111}
    \InternalNode[x=8]{11111111}
    
    \ExternalNode[x=0,y=-1,label={0,1},position=below]{0}
    \ExternalNode[x=1,y=-1,label={2,3},position=below]{10}
    \ExternalNode[x=2,y=-1,label={4,5},position=below]{110}
    \ExternalNode[x=3,y=-1,label={6,7},position=below]{1110}
    \ExternalNode[x=4,y=-1,label={8,9},position=below]{11110}
    \ExternalNode[x=5,y=-1,label={10,11},position=below]{111110}
    \ExternalNode[x=6,y=-1,label={12,13},position=below]{1111110}
    \ExternalNode[x=7,y=-1,label={14,15},position=below]{11111110}
    \ExternalNode[x=8,y=-1,label={16,17},position=below]{111111110}
    \ExternalNode[x=9,y=0,label={18,19},position=below]{111111111}

    \Edge(D)(root)
    \Edge(root)(0)
    \Edge(root)(1)
    \Edge(1)(10)
    \Edge(1)(11)
    \Edge(11)(110)
    \Edge(11)(111)
    \Edge(111)(1110)
    \Edge(111)(1111)
    \Edge(1111)(11110)
    \Edge(1111)(11111)
    \Edge(11111)(111110)
    \Edge(11111)(111111)
    \Edge(111111)(1111110)
    \Edge(111111)(1111111)
    \Edge(1111111)(11111110)
    \Edge(1111111)(11111111)
    \Edge(11111111)(111111110)
    \Edge(11111111)(111111111)
\end{tikzpicture}
        \caption{Tensor train partition (TT)\label{fig:cascade-partition-TT}}
    \end{subfigure}
    \\[4ex]
    \begin{subfigure}[b]{\textwidth}
        \centering
        \begin{tikzpicture}[rotate=180]        
    \InternalNode{root}
    \node[above=0.5cm of root](D){};

    \InternalNode[x=-3,y=1]{1}
    \InternalNode[x=3,y=1]{0}
    
    \InternalNode[x=-4.5,y=2]{11}
    \InternalNode[x=-1.5,y=2]{10}
    \InternalNode[x=1.5,y=2]{01}
    \InternalNode[x=4.5,y=2]{00}

    \InternalNode[x=-5.25,y=3]{111}
    \ExternalNode[x=-3.75,y=3,label={14,15},position=below]{110}
    \ExternalNode[x=-2.25,y=3,label={12,13},position=below]{101}
    \ExternalNode[x=-0.75,y=3,label={10,11},position=below]{100}
    \ExternalNode[x=0.75,y=3,label={8,9},position=below]{011}
    \ExternalNode[x=2.25,y=3,label={6,7},position=below]{010}
    \ExternalNode[x=3.75,y=3,label={4,5},position=below]{001}
    \InternalNode[x=5.25,y=3]{000}

    \ExternalNode[x=-6.0,y=4,label={18,19},position=below]{1111}
    \ExternalNode[x=-4.5,y=4,label={16,17},position=below]{1110}
    \ExternalNode[x=4.5,y=4,label={2,3},position=below]{0001}
    \ExternalNode[x=6.0,y=4,label={0,1},position=below]{0000}

    \Edge(D)(root)

    \Edge(root)(0)
    \Edge(root)(1)

    \Edge(0)(00)
    \Edge(0)(01)
    \Edge(1)(10)
    \Edge(1)(11)

    \Edge(00)(000)
    \Edge(00)(001)
    \Edge(01)(010)
    \Edge(01)(011)
    \Edge(10)(100)
    \Edge(10)(101)
    \Edge(11)(110)
    \Edge(11)(111)

    \Edge(000)(0000)
    \Edge(000)(0001)
    \Edge(111)(1110)
    \Edge(111)(1111)
\end{tikzpicture}
        \caption{Binary tree partition (BT)\label{fig:cascade-partition-BT}}
    \end{subfigure}
    \caption{Tree structures of the partitions of the reaction cascade example. For both partitions 9 pathways are severed. Small numbers at the leaves of the tree structure denote on which species the corresponding low-rank factors depend.\label{fig:cascade-partitions}}
\end{figure}
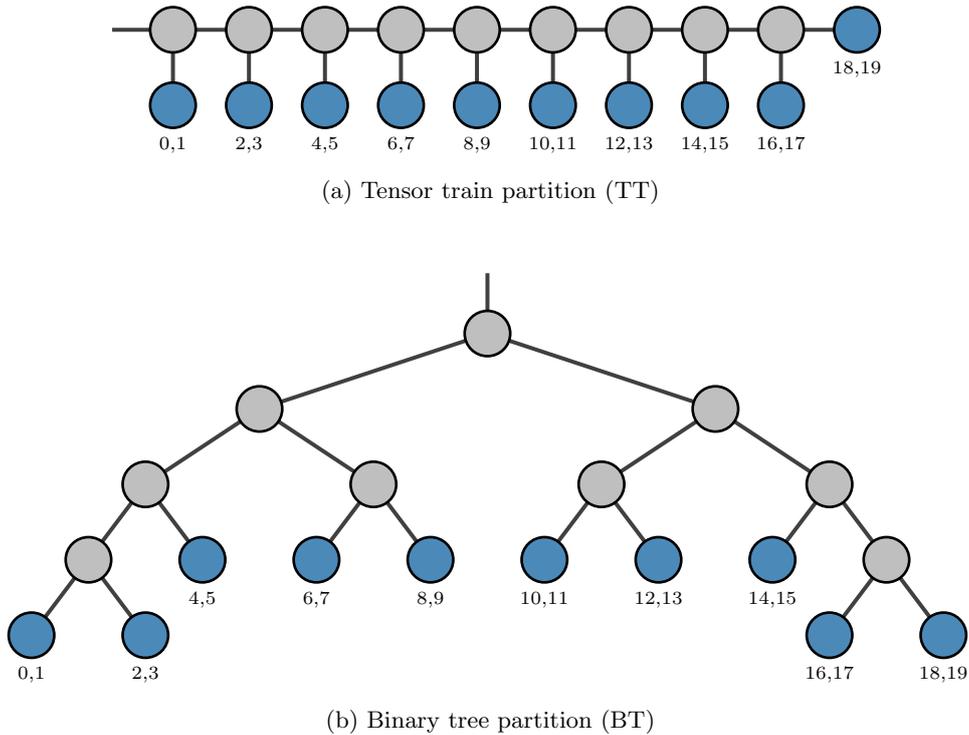

To compare our results with SSA, we compute solutions with the PS-TTN integrator using the tensor train partition (TT) shown in Figure~\ref{fig:cascade-partition-TT}. Since we want to severe only a few pathways of the reaction network while at the same time minimizing the number of degrees of freedom of the low-rank factors, we always keep two adjacent species in a single leaf. Thereby a single-low rank factor has $64^2=4096$ degrees of freedom and in total only 9 pathways are severed. This approach also helps to keep the height of the tree relatively small, which results in small recursion depths for the algorithm. For the time integration an explicit Euler scheme with time step size $\Delta t=10^{-1}$ was employed. Due to the extremely large system size with $64^{20}=1.3 \cdot 10^{36}$ degrees of freedom the full CME can no longer be solved, we therefore compare the PS-TTN solutions with a reference solution obtained via SSA using $10^7$ runs. All computations were performed on the same computer system as in the previous example.

\begin{figure}[H]
    \centering
    \includegraphics[scale=0.8]{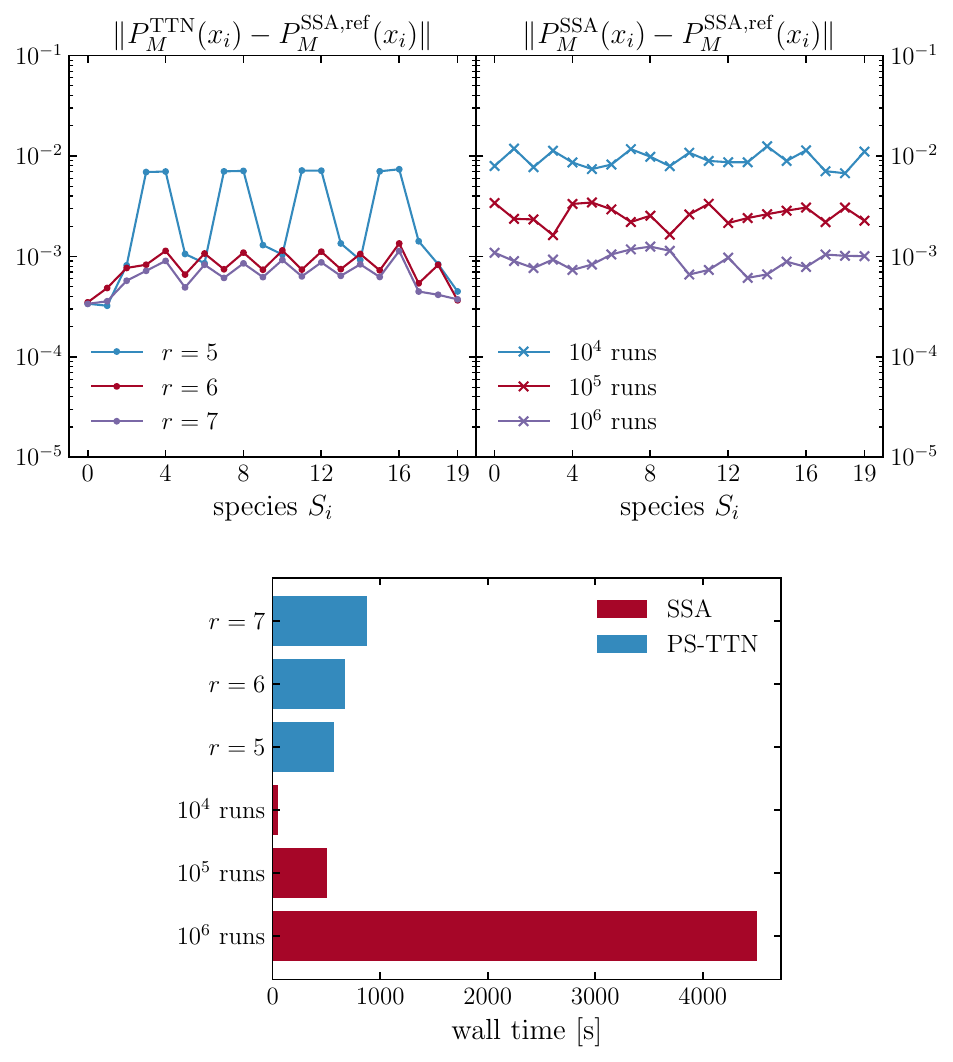}
    \caption{Top:~2-norm error of the marginal distributions for the reaction cascade model using the PS-TTN integrator (left panel) and SSA (right panel). For the PS-TTN integrator the TT partition of Figure~\ref{fig:cascade-partition-TT} was chosen. The reference solution was obtained with SSA using $10^7$ runs. For the PS-TTN we used explicit Euler with time step size $\Delta t=10^{-1}$. Bottom:~Wall time comparison (in seconds) for the PS-TTN integrator and SSA.\label{fig:cascade-ttn-ssa}}
\end{figure}

In the top panel of Figure~\ref{fig:cascade-ttn-ssa}, we again compare the accuracy of our method with SSA. Due to the large system size we now calculate the 2-norm error of all marginal distributions instead of the full probability distribution. For the PS-TTN integrator we use equal ranks for all nodes.

We observe that for both rank $r=6$ and $r=7$ the error of the PS-TTN solution is similar as for the Monte Carlo method with $10^6$ runs, but with the advantage of being approximately four times faster. In other words, our method is as fast as SSA with $10^5$ runs, but shows a similar accuracy as SSA with $10^6$ runs. Recall that the reference solution is not exact, therefore the real error for the PS-TTN solution with $r=7$ might even be smaller. Also note that although SSA with $10^4$ is one order of magnitude faster than the PS-TTN integrator, the solution is very inaccurate and due to the $1/\sqrt{N}$ scaling of Monte Carlo methods $10^6$ runs yield only an improvement of a factor of $10$.

As an aside, the error of the PS-TTN solution with rank $r=5$ exhibits interesting periodic features: After the first three points it alternates between larger values of approximately $10^{-2}$ and smaller ones of $10^{-3}$. We believe that the periodicity of this behaviour is related to the fact that all low-rank factors on the leaf nodes depend on two species. It seems reasonable that the higher errors occur always for the two species which are in different partitions and where the reaction is approximated.

The memory requirements for storing the full probability distribution at a single point in time are $1.06 \cdot 10^{31}\,\mathrm{MB}$, whereas for the low-rank approximation with the TTN partition only $2.3\,\mathrm{MB}$ for $r=7$ are needed. This impressively demonstrates the compression capabilities of tree tensor networks.

\begin{figure}[H]
    \centering
    \includegraphics[scale=0.8]{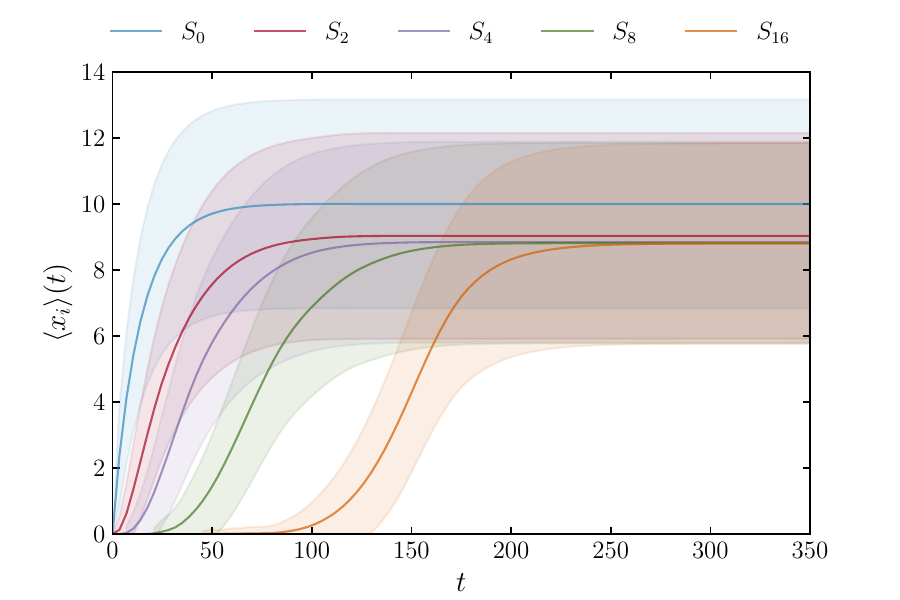}
    \caption{Time-dependent mean population numbers with standard deviation (shaded areas) for five species of the reaction cascade model.\label{fig:cascade-concentrations}}
\end{figure}

In some applications primarily time-dependent mean population numbers $\langle x_i \rangle(t)$ ($i=0,\dots,d-1$) (or concentrations, when dividing by the volume of the system) are of interest. The mean population numbers are first moments of the probability distribution,
\begin{equation*}
    \langle x_i \rangle (t) = \sum_{\vb{x} \in \Omega_{\vb{\zeta},\vb{\eta}}} x_i P(t, \vb{x}),
\end{equation*}
and they can be efficiently computed from the tree tensor representation of the probability distribution. In Figure~\ref{fig:cascade-concentrations} we plot the time-dependent mean population numbers along with the standard deviation for five different species. Since the propensity functions for the creation of species $S_i$ depends on the population number of $S_{i-1}$, the growth of the population numbers is delayed in time. The standard deviation also increases over time and from time $t=200$ it is very high for all species. This implies that stochastic fluctuations are non-negligible and therefore a deterministic description of this system is not feasible.

\begin{figure}[H]
    \centering
    \includegraphics[scale=0.8]{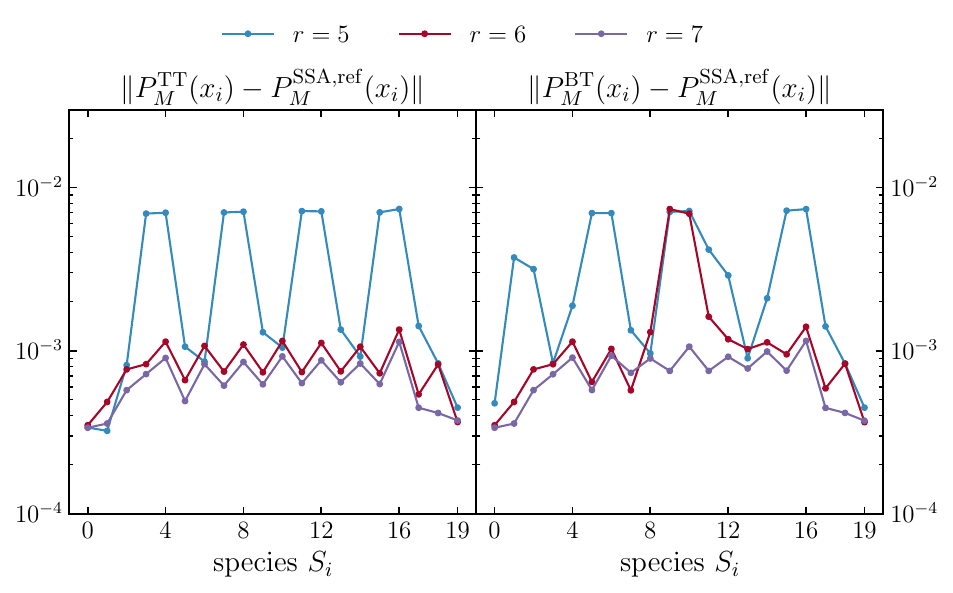}
    \caption{2-norm error of the marginal distributions for the reaction cascade model using the PS-TTN integrator with the TT (left panel) and BT partition (right panel). The reference solution was obtained with SSA using $10^7$ runs. For both partitions we used the explicit Euler scheme with time step size $\Delta t=10^{-1}$.\label{fig:cascade-partition-comparison}}
\end{figure}

Finally we compare the TT partition with the binary tree partition shown in Figure~\ref{fig:cascade-partition-BT}. We again use the explicit Euler scheme with time step size $\Delta t = 10^{-1}$ and use as a reference solution SSA with $10^7$ runs. From Figure~\ref{fig:cascade-partition-comparison} we conclude that for $r=5$ both solutions show approximately the same accuracy. For $r=6$ the TT solution is more accurate, whereas for $r=7$ both solutions have a similar accuracy. Overall, the TT partition seems to be more suitable for this particular example.

\section{Conclusion and outlook}\label{sec:outlook}
In this work we presented an integrator for the kinetic chemical master equation based on the projector-splitting integrator for tree tensor networks of \cite{Ceruti_2020}. We implemented the numerical scheme efficiently by exploiting the factorization property and the reactant dependency of the propensity functions. The resulting integrator overcomes the curse of dimensionality and thus can outperform state-of-the-art Monte Carlo techniques in terms of run time and accuracy even for large problems such as the 20-dimensional reaction cascade of Section~\ref{sec:cascade}, with the added advantage that our method is completely noise-free.

Our implementation can be easily extended to a rank-adaptive scheme for TTNs, for example to the rank-adaptive BUG integrator of \cite{Ceruti_2023}. Furthermore, the sliding windows technique of \cite{Henzinger_2009,Wolf_2010} could be implemented, where the truncated state space is adapted according to the time evolution of the probability distribution. Automatizing the partitioning of the reaction network is a subject of ongoing research, and in future we also want to use our integrator to study large biochemical reaction networks and reaction-diffusion processes.

\section{Acknowledgements}
The authors thank Gianluca Ceruti and Dominik Sulz for fruitful discussions.

\printbibliography

\end{document}